\documentclass[10pt]{amsart}

\usepackage{verbatim}
\usepackage{amssymb}
\usepackage{graphicx, subfigure}

\begin{document}


\newcommand{\bR}{\mathbb{R}}	


\newcommand{\A}{\alpha}
\newcommand{\B}{\beta}
\newcommand{\G}{\gamma}





\newtheorem{main}{Theorem}
\renewcommand{\themain}{\Alph{main}}
\newtheorem{thm}{Theorem}[section]
\newtheorem{remark}{Remark}[section]
\newtheorem{lemma}[thm]{Lemma}
\newtheorem{prop}[thm]{Proposition}






\title[Apollonian Sets]{Apollonian sets in taxicab geometry}
\author[SUMmER REU 2017]{Eric Bahuaud,
Shana Crawford,
Aaron Fish,
Dylan Helliwell,
Anna Miller,
Freddy Nungaray,
Suki Shergill,
Julian Tiffay,
Nico Velez
}
\address{Department of Mathematics, 
Seattle University.}
\email{bahuaude(at)seattleu.edu}
\email{helliwed(at)seattleu.edu}

\date{\today}

\begin{abstract}
Fix two points $p$ and $q$ in the plane and a positive number $k \neq 1$. A result credited to Apollonius of Perga states that the set of points $x$ that satisfy $d(x,p)/d(x,q) = k$ forms a circle.  In this paper we study the analogous set in taxicab geometry.  We find that while Apollonian sets are not taxicab circles, more complicated Apollonian sets can be characterized in terms of simpler ones.
\end{abstract}

\maketitle


\section{Introduction}

The taxicab plane is the set $\bR^2$ endowed with the $\ell^1$ or taxicab metric given by
\[
d(x, y) = |x_1 - y_1| + |x_2 - y_2|.
\]

There is by now a long list of papers that study the differences between particular notions in Euclidean geometry and their taxicab counterparts.  For a non-exhaustive list, see for example \cite{Krause}, where taxicab geometry was introduced as a tool for training in research mathematics; \cite{Thompson} where comparisons between area and angles are made; \cite{Reynolds} and \cite{KAGO} for two examples of the exploration of taxicab conics; and \cite{Kaya} for comparisons with other classical Euclidean theorems.

\begin{figure}
\begin{picture}(360,120)
\put(5,20){
\includegraphics[scale = .3, clip = true, draft = false]{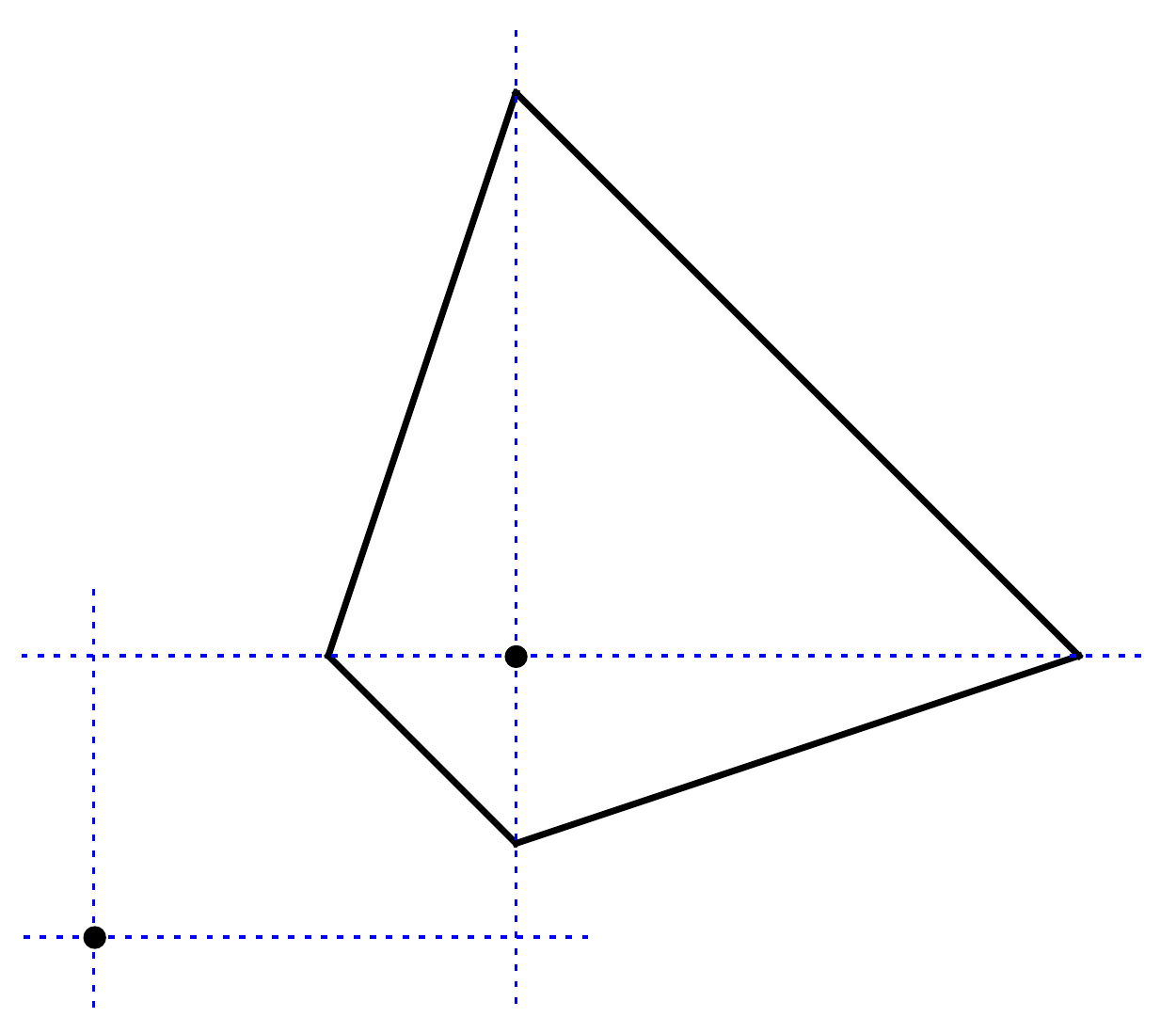}
}
\put(19,21){$p$}
\put(58,48){$q$}
\put(50,0){(a)}

\put(125,20){
\includegraphics[scale = .3, clip = false, draft = false]{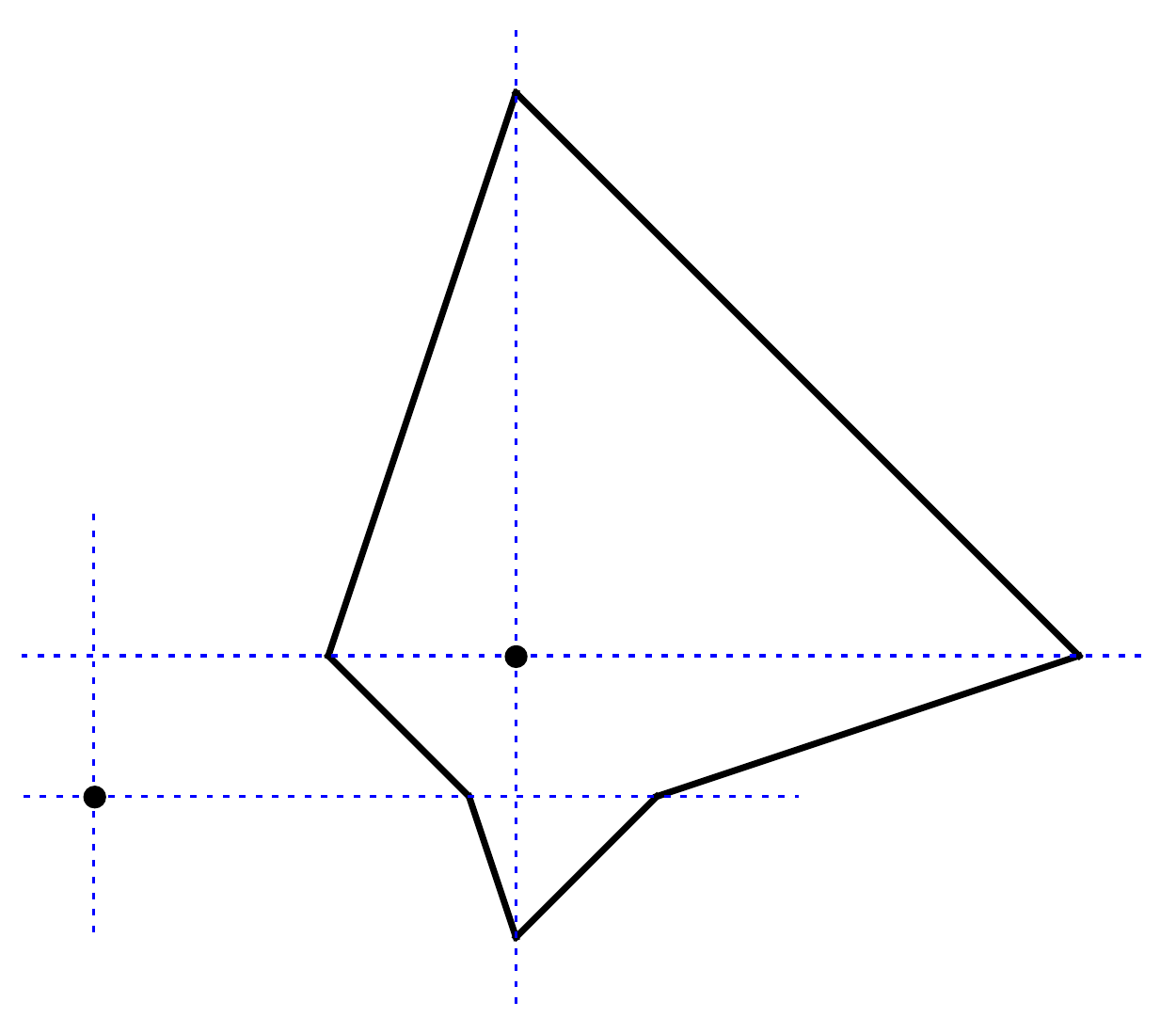}
}
\put(139,35){$p$}
\put(178,48){$q$}
\put(170,0){(b)}

\put(245,10){
\includegraphics[scale = .3, clip = false, draft = false]{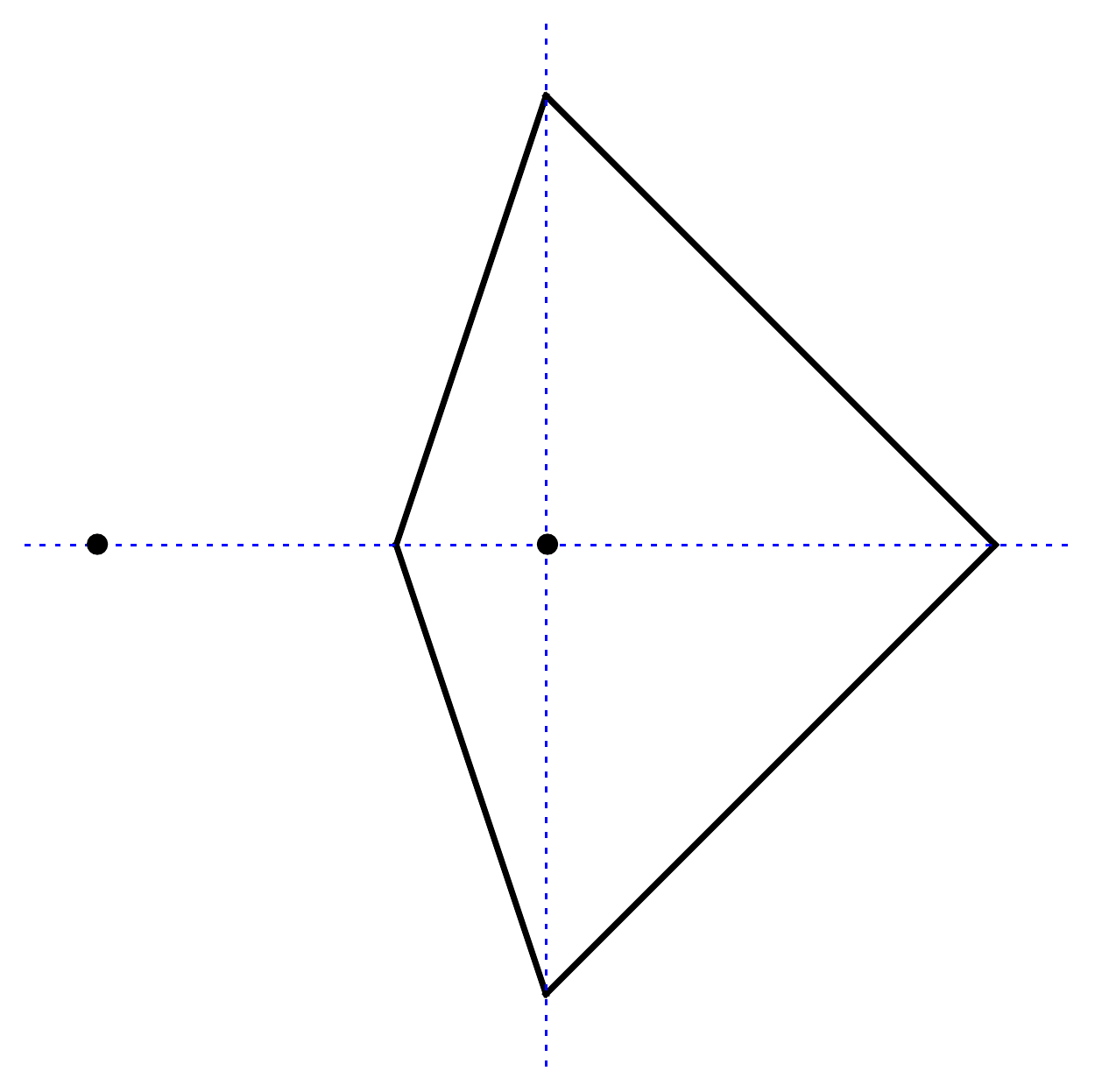}
}
\put(258,56){$p$}
\put(305,56){$q$}
\put(290,0){(c)}

\end{picture}

\caption{Typical Apollonian sets $A(p, q; k)$.  In all three cases, $p = (0,0)$ and $k = 2$: (a) $q = (3,2)$, (b) $q = (3,1)$, (c) $q = (3,0)$} \label{examplefig}
\end{figure}
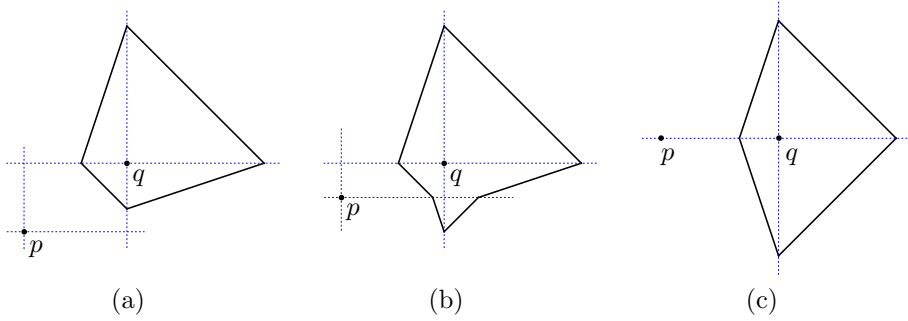

In this paper we are interested in a classical construction in Euclidean geometry, attributed to the 3rd century BCE mathematician Apollonius of Perga.  Fix two distinct points $p, q$ in the plane and a positive real constant $k$.  Consider the set
\[
A(p,q;k)= \left\{ x \in \bR^2: \frac{d(x,p)}{d(x,q)} = k \right\}.
\]
If $d$ is the standard Euclidean metric, then Apollonius's result states that $A(p,q;k)$ forms either a circle if $k \neq 1$, or a straight line if $k=1$.  The definition of $A(p,q;k)$ makes sense in any metric space $(X,d)$ and in this paper we characterize the Apollonian sets $A(p,q;k)$ for $\bR^2$ endowed with the taxicab metric.  We find that in no case does the set coincide with a (taxicab) circle, but instead takes on a number of different shapes depending on the relative positions of $p$ and $q$ and the value of $k$.  See Figure \ref{examplefig}.  Despite this apparent complexity, we find that, with the exception of $k = 1$, all Apollonian sets can be expressed as a union of trapezoids that are themselves Apollonian sets.  A considerable part of this paper is dedicated to introducing notation that we hope will further the development of taxicab geometry from a more synthetic viewpoint.  To date many of the papers on taxicab geometry exploit the piecewise linear structure of the taxicab metric, leading to the necessity of detailed case analysis.  Our definitions and lemmas are designed so that few such algebraic computations are required.

We now state our results more precisely, referring the reader to the next section for more background.  First, in order to study Apollonian sets, it is valuable to think of $A(p,q;k)$ as the boundary of a compact closed region which, for $k > 1$, has the form
\[
B(p,q;k) = \left\{ x \in \bR^2: \frac{d(x,p)}{d(x,q)} \geq k \right\}.
\]
Second, in view of the isometries of the taxicab plane, horizontal and vertical lines which we term \textit{coordinate lines}, and lines with slope $\pm 1$ which we term \textit{guide lines} play an important role.  In particular given $p$ and $q$ in $\bR^2$, the guide lines emanating from $p$ and $q$ respectively form a rectangle with vertices $p$, $q$ and two other points called \textit{guide complements} that we label $g^+$ and $g^-$.  We prove for $k \neq 1$ that the filled Apollonian sets $B(g^+,q;k)$ and $B(g^-,q;k)$ are both filled isosceles trapezoids and that other filled Apollonian sets are unions of these.  Our two main results are as follows:

\begin{main} \label{trapthm}
Let $p$ and $q$ share a guide line $gl$ and, without loss of generality, let $k > 1$.  Then the Apollonian set $A(p, q; k)$ is an isosceles trapezoid with the following properties:
\begin{itemize}
\item the line of symmetry for $A(p, q; k)$ is $gl$;
\item the vertices of $A(p, q; k)$ all lie on the coordinate lines of $q$;
\item if the legs of $A(p, q; k)$ are extended to lines, these lines intersect at $p$;
\item the slopes of the legs are $m \frac{k+1}{k-1}$ and $m \frac{k-1}{k+1}$, where $m$ is the slope of $gl$.
\end{itemize}
\end{main}

\begin{main} \label{apollonianisunionthm}
Let $p,\ q \in \bR^2$, let $g^+$ and $g^-$ be the guide complements of $p$ and $q$, and let $k \in [0, 1) \cup (1, \infty]$.  Then
\[ B(p, q; k) = B(g^+, q; k) \cup B (g^-, q; k). \]
\end{main}

Since $B(g^{\pm},q;k)$ are characterized by Theorem \ref{trapthm}, these two results completely characterize Apollonian sets when $k \neq 1$ and provide a constructive method for producing them.  The sets that arise when $k = 1$ have been studied in different contexts, see \cite{Reynolds} and \cite{KAGO}, and are also included here for completeness.

This paper is structured as follows: in Section \ref{sec:back} we introduce taxicab geometry, our notation, various reference objects, and other background material including a brief discussion of isometries.  In Section \ref{sec:as} we introduce Apollonian sets in the taxicab plane, recognize that these sets fit into a more general framework, and establish a number of foundational facts about them.  Finally, in Section \ref{sec:char}, we provide a complete characterization of all Apollonian sets, including proofs of the two main theorems above.

This paper reports on some of the work done in the Seattle University SUMmER REU, 2017.  We gratefully acknowledge the support of Seattle University's NSF REU \#1460537 (SUMmER in Seattle).  Additionally, the work of A. Miller and A. Fish were supported by a generous gift to Seattle University from Ms. Rose Southall.

\section{Background} \label{sec:back}

In this section we introduce the distance function and give some basic background in taxicab geometry.  The reader may also consult \cite{Krause}.  We also introduce several special sets of interest in the paper and discuss taxicab isometries.

\subsection{Distance Functions}

A distance function, or metric, on a set $A$ is a function
\[
d: A \times A \longrightarrow \mathbb{R}
\]
satisfying
\begin{enumerate}
\item for all  $x, y \in A$, $d(x, y) \geq 0$ and $d(x, y) = 0$ if and only if $x = y$;
\item for all $x, y \in A$, $d(x, y) = d(y, x)$;
\item for all $x, y, z \in A$, $d(x,z) \leq d(x, y) + d(y, z)$.
\end{enumerate}

A set together with a distance function is called a metric space.

In this paper, the set $A$ will be the plane $\mathbb{R}^2$ and we denote a point  $p \in \bR^2$ in components by $p = (p_1, p_2)$.  The most familiar distance function on $\bR^2$ is the Euclidean distance function $d_E$, given by
\[
d_E(x, y) = \bigl((x_1 - y_1)^2 + (x_2 - y_2)^2\bigr)^{\frac{1}{2}}.
\]
Alternatively, the taxicab distance function $d$ is given by
\[
d(x, y) = |x_1 - y_1| + |x_2 - y_2|.
\]
Unless otherwise specified, an unadorned $d$ will be the taxicab distance.

\subsection{Special points, lines, regions}

As we will see in the remainder of the paper, it will be helpful to introduce a few special reference objects that occur frequently when discussing taxicab geometry.

A fundamental object in taxicab geometry is the circle.  The circle centered at $p$ and with radius $r > 0$ is a square with vertices at $(p_1 \pm r, p_2)$ and $(p_1, p_2 \pm r)$.  It will turn out in this paper that circles do not make a significant appearance, and it is especially worth noting that Apollonian sets turn out never to be circles.

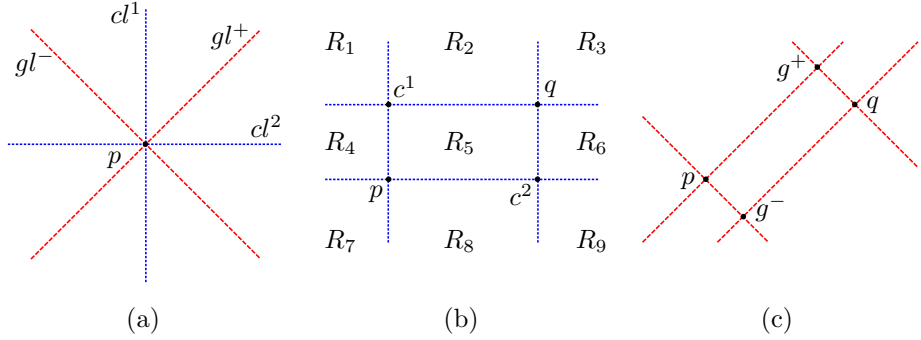
\begin{figure}
\begin{picture}(360,120)
\put(0,10){
\includegraphics[scale = .3, draft = false]{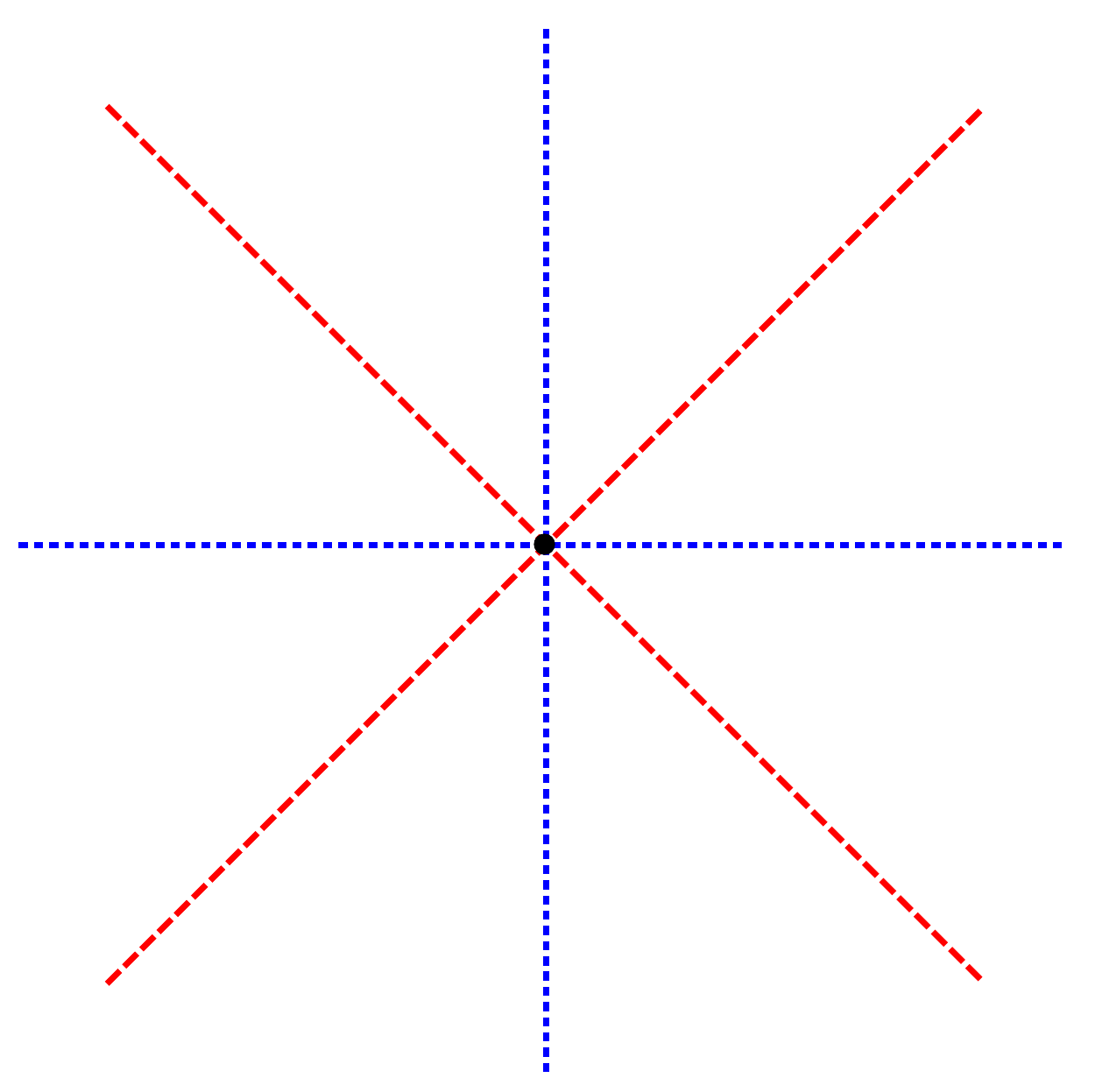}
}
\put(43,57){$p$}
\put(97,67){$cl^2$}
\put(44,110){$cl^1$}
\put(82,102){$gl^+$}
\put(8,92){$gl^-$}
\put(50,-5){(a)}

\put(120,25){
\includegraphics[scale = .3, draft = false]{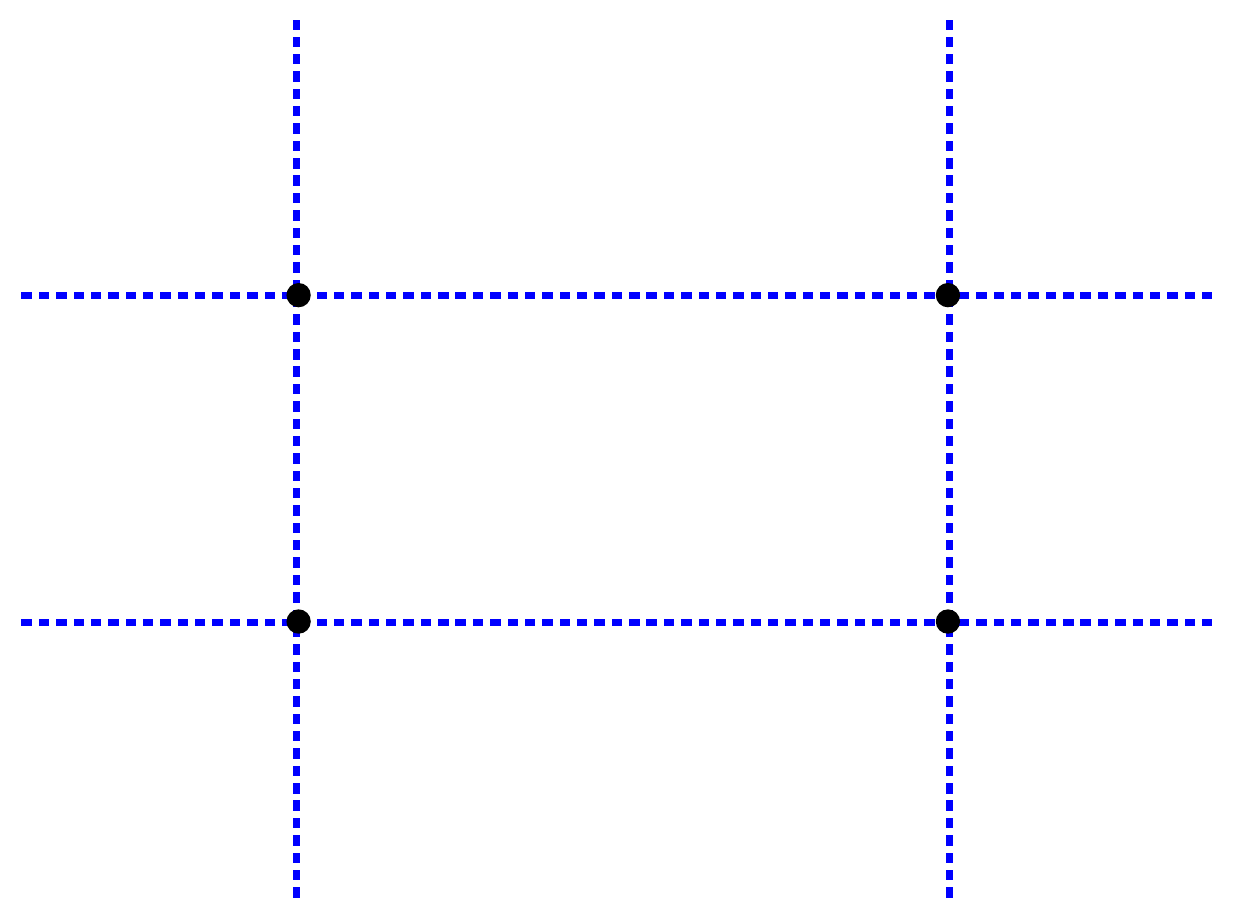}
}
\put(142,44){$p$}
\put(208,84){$q$}
\put(151,82){$c^1$}
\put(195,41){$c^2$}
\put(125, 100){$R_1$}
\put(170, 100){$R_2$}
\put(220, 100){$R_3$}
\put(125, 62){$R_4$}
\put(170, 62){$R_5$}
\put(220, 62){$R_6$}
\put(125, 25){$R_7$}
\put(170, 25){$R_8$}
\put(220, 25){$R_9$}
\put(170,-5){(b)}

\put(240,25){
\includegraphics[scale = .3, draft = false]{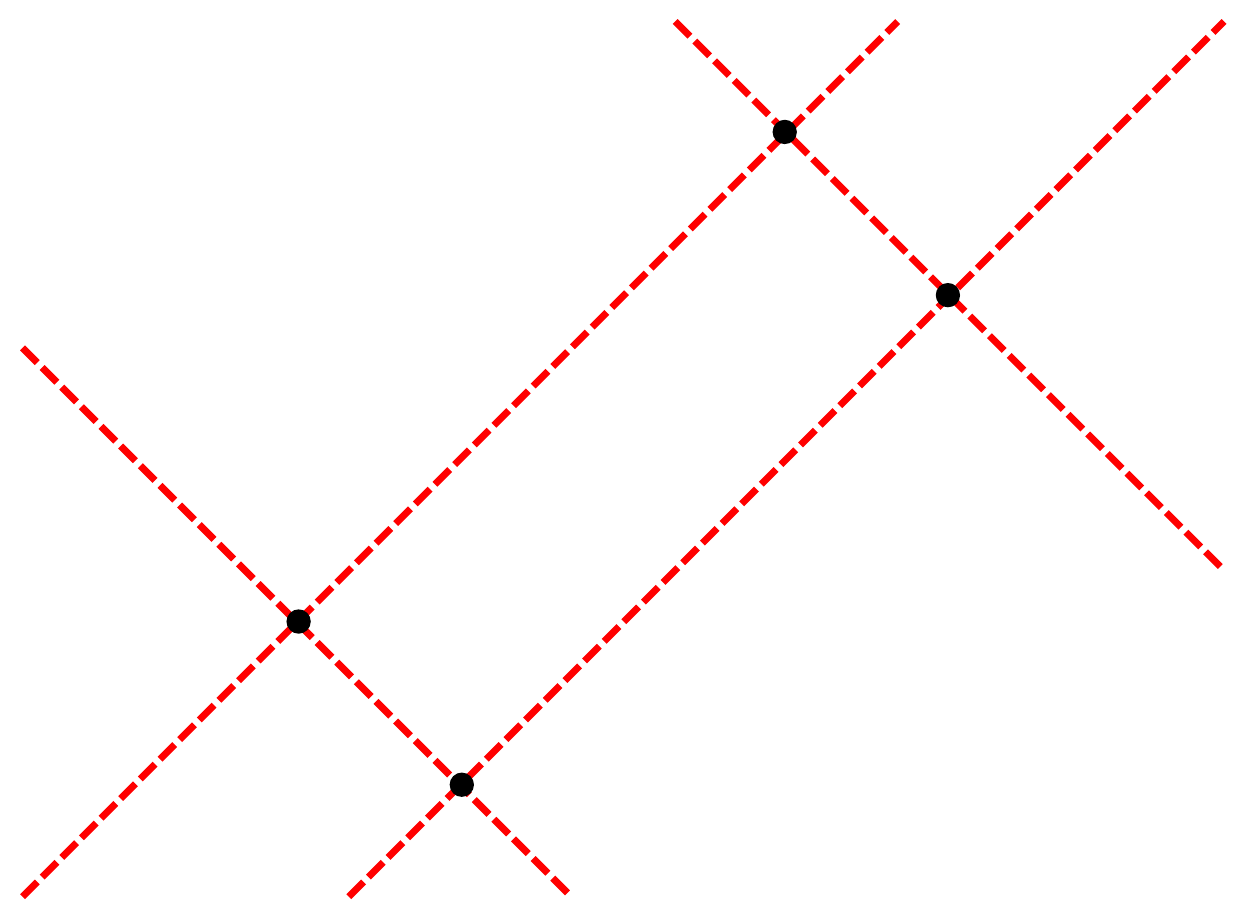}
}
\put(260,49){$p$}
\put(330,78){$q$}
\put(296,90){$g^+$}
\put(288,37){$g^-$}
\put(290,-5){(c)}

\end{picture}

\caption{Reference objects:  (a) coordinate lines and guide lines of $p$, (b) coordinate complements ($c^i = c^i(p, q)$) and regions, (c) guide complements ($g^{\pm} = g^{\pm}(p, q)$)} \label{linefig}
\end{figure}

Given a point $p = (p_1, p_2)$, we define the \textit{coordinate lines} of $p$ to be the lines
\begin{align*}
cl^1(p) &= \{(x_1, x_2): x_1 = p_1\} = \{(p_1, x_2)\} \\
cl^2(p) &= \{(x_1, x_2): x_2 = p_2\} = \{(x_1, p_2)\}.
\end{align*}
Note that if $p = (0, 0)$, then the coordinate lines are just the usual coordinate axes.  See Figure \ref{linefig}(a).  Given a point $q$, the coordinate lines of $q$ divide $\bR^2$ into four quadrants we call $q$-quadrants.  Unless otherwise mentioned, these quadrants are to be interpreted as including their boundaries.

We define the \textit{coordinate complements} of $p$ and $q$ to be the points
\begin{align*}
c^1(p,q) &= cl^1(p) \cap cl^2(q) = (p_1, q_2) \\
c^2(p,q) &= cl^2(p) \cap cl^1(q) = (q_1, p_2).
\end{align*}
Together, $p$, $q$, $c^1(p,q)$, and $c^2(p,q)$ form the vertices of a filled rectangle whose edges lie on the coordinate lines of $p$ and $q$.  We call this filled rectangle the \textit{coordinate rectangle}, and we note that this rectangle collapses to a line segment if $p$ and $q$ share a coordinate line.  The coordinate rectangle plays an analogous role to a line segment in Euclidean geometry with respect to the triangle inequality: every point $z$ in the coordinate rectangle satisfies the equality
\[ d(p,q) = d(p,z) + d(z,q). \]
The coordinate lines of $p$ and $q$ divide $\bR^2$ into 9 regions $R_i,\ i = 1, \ldots 9$, as indicated in Figure \ref{linefig}(b).  These regions are always labelled as indicated, independent of the relative positions of $p$ and $q$.  Region 5 is the coordinate rectangle.  As with quadrants, these regions are to be interpreted as including their boundaries.

Next, we define the \textit{guide lines} of $p$ to be the lines
\begin{align*}
gl^+(p) &= \{(x_1, x_2): x_2 - p_2 = x_1 - p_1\} \\
gl^-(p) &= \{(x_1, x_2): x_2 - p_2 = -(x_1 - p_1)\}.
\end{align*}
These are the lines passing through $p$ with slope $\pm 1$.  See Figure \ref{linefig}(a).

We also define the \textit{guide complements} of $p$ and $q$ to be the points
\begin{align*}
g^+(p,q) &= gl^+(p) \cap gl^-(q) \\
g^-(p,q) &= gl^-(p) \cap gl^+(q).
\end{align*}
Together, $p$, $q$, $g^+(p,q)$, and $g^-(p,q)$ form the vertices of a rectangle whose edges lie on the guide lines of $p$ and $q$.  See Figure \ref{linefig}(c).  For both the coordinate complements and the guide complements, the superscripts have been chosen to reflect the corresponding line determined by the first argument.

Guide complements play a significant role in understanding Apollonian sets.  The following two lemmas establish important relationships between guide complements and the points defining them.

\begin{lemma} If $p$ and $q$ do not share a guide line, then the guide complements $g^+(p, q)$ and  $g^-(p, q)$ lie outside the coordinate rectangle of $p$ and $q$.  If $p$ and $q$ do share a guide line, then guide complements are $p$ and $q$.
\end{lemma}

The second statement in this lemma is analogous to the fact that if $p$ and $q$ share a coordinate line, then the coordinate complements are $p$ and $q$.

\begin{proof}
Suppose $p$ and $q$ do not share a guide line.  One of the guide lines of $p$ will intersect the coordinate rectangle only at $p$.  The guide complement defined by this line must then lie outside the rectangle.  To see that both guide complements lie outside, note that the guide line of $q$ that is parallel to the rectangle-avoiding guide line of $p$ must also avoid the rectangle, and this line is involved in defining the second guide complement.

The second statement in the lemma follows immediately from the definition of guide complement.
\end{proof}

\begin{lemma} \label{guidequadlemma}
Let $g$ be either of the guide complements of $p$ and $q$.  If $p$ and $q$ do not share a guide line, then $p$ and $q$ lie in different $g$-quadrants.
\end{lemma}

\begin{proof}
Suppose not.  Then $p$ and $q$ are in the same $g$-quadrant.  But this implies that the slope of the line formed by $g$ and $p$ has the same sign as the slope of the line formed by $g$ and $q$, a contradiction.
\end{proof}

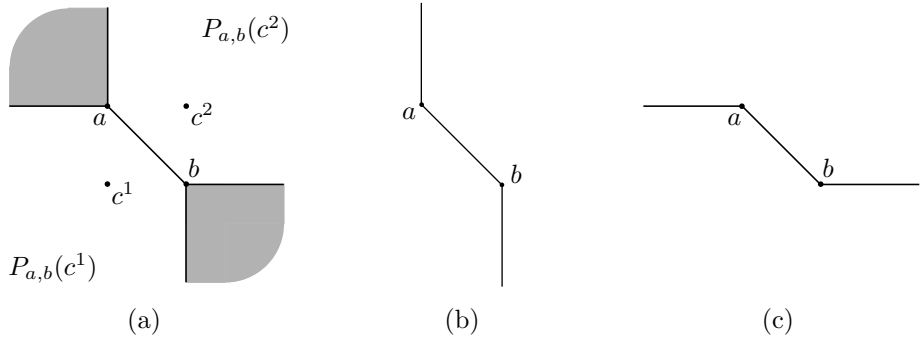
\begin{figure}
\begin{picture}(360,120)
\put(0,10){
\includegraphics[scale = .3, clip = true, draft = false]{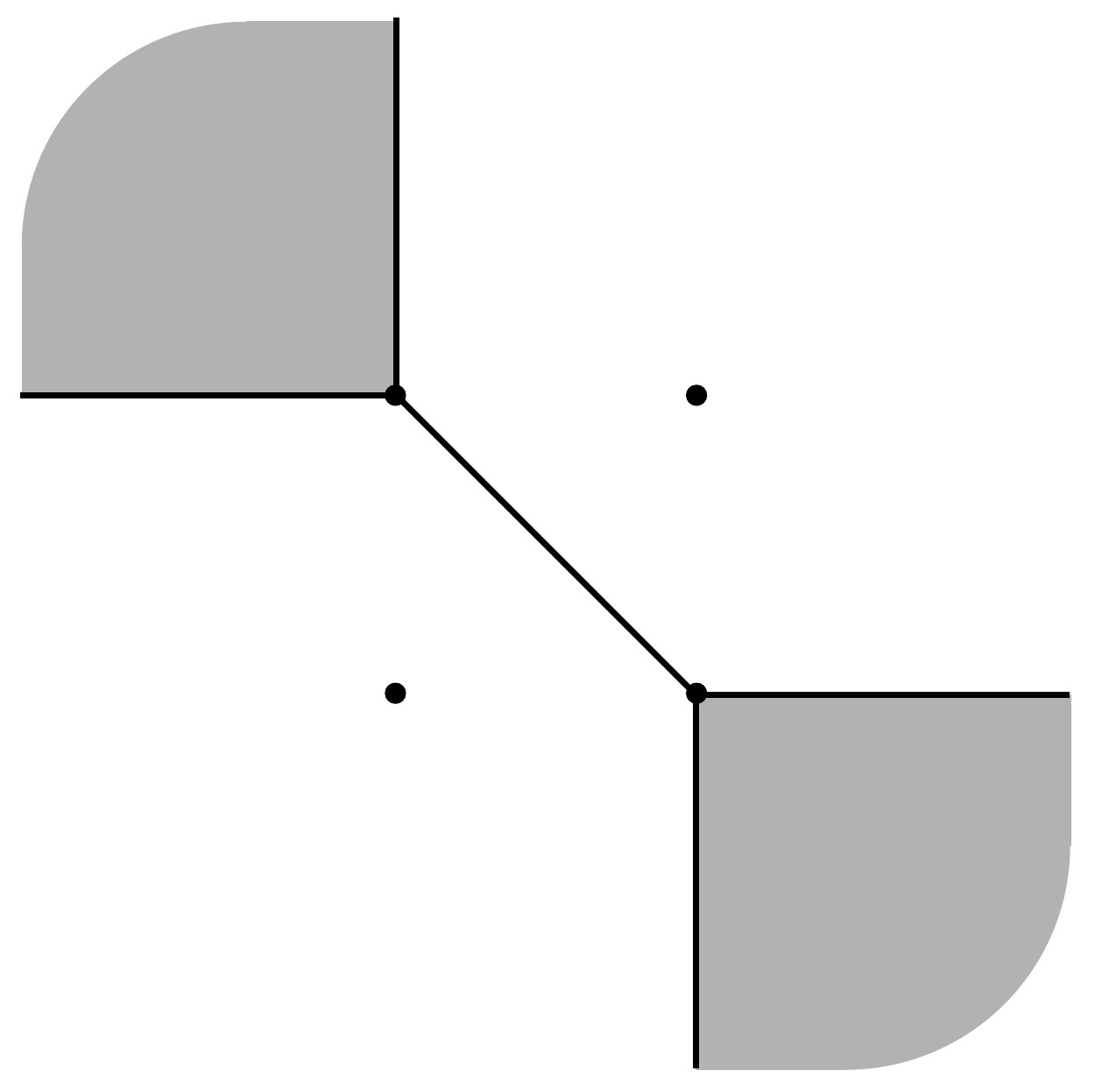}
}
\put(37,72){$a$}
\put(73,52){$b$}
\put(74,72){$c^2$}
\put(44,42){$c^1$}
\put(78,105){$P_{a,b}(c^2)$}
\put(5,15){$P_{a,b}(c^1)$}
\put(50,-5){(a)}

\put(157,10){
\includegraphics[scale = .25, clip = false, draft = false]{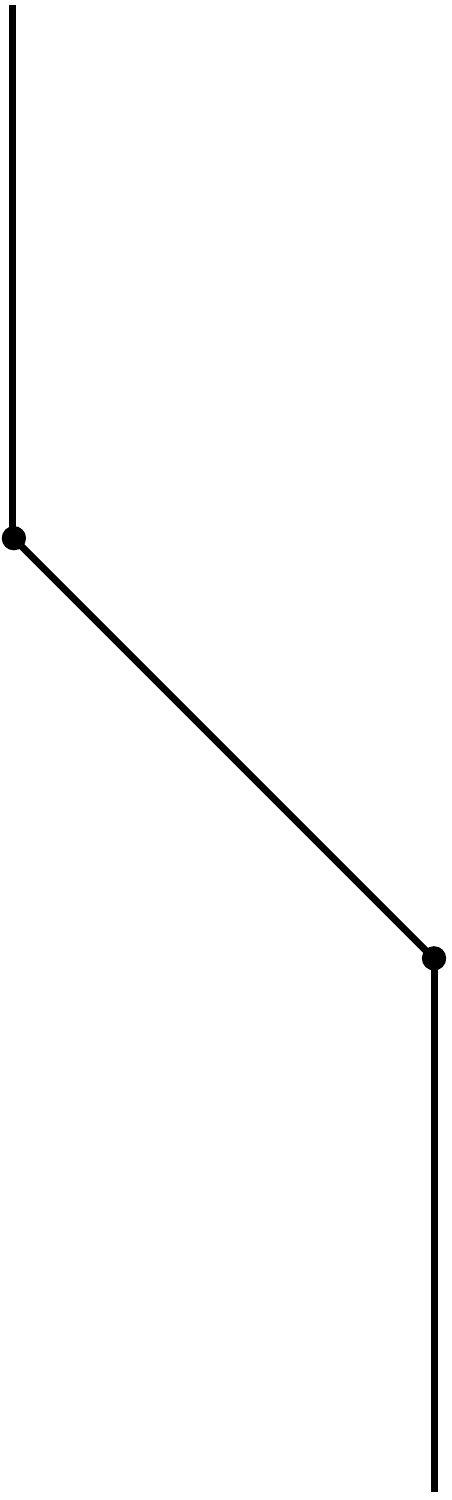}
}
\put(154,74){$a$}
\put(195,50){$b$}
\put(170,-5){(b)}

\put(240,47){
\includegraphics[scale = .3, clip = false, draft = false]{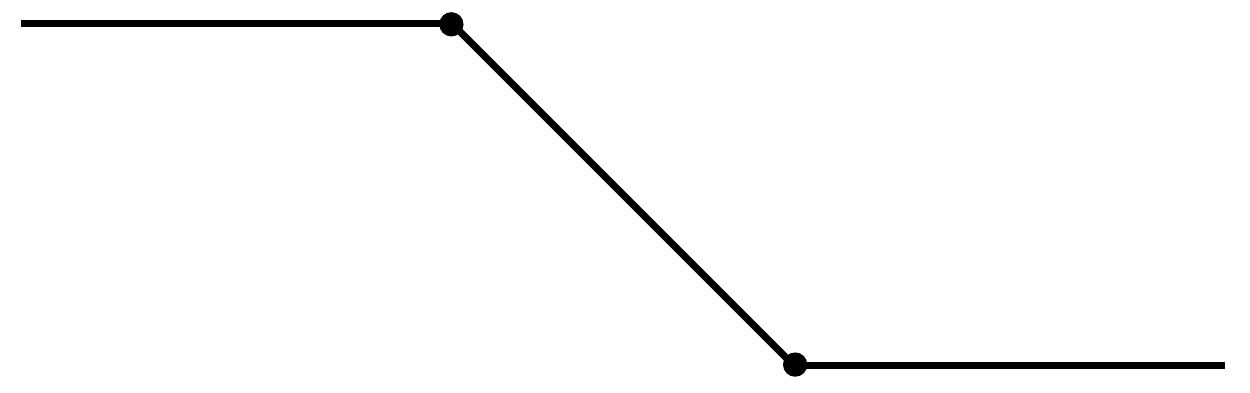}
}
\put(277,72){$a$}
\put(313,52){$b$}
\put(290,-5){(c)}

\end{picture}

\caption{Barbells and lightning bolts:  (a) barbell of $a$ and $b$, and $P_{a,b}(c^i)$ (unshaded components); (b) lightning bolt of type 1; (c) lightning bolt of type 2} \label{barbellfig}
\end{figure}

Let $a$ and $b$ share a guide line $gl$.  Then the \textit{barbell} $bb(a, b)$ is the union of the $a$-quadrant that contains the half of $gl$ that does not include $b$, the $b$-quadrant that contains the half of $gl$ that does not include $a$, and $gl$.  The complement of $bb(a, b)$ comprises two components, each of which contains one of the coordinate complements $c^i$, $i = 1, 2$, of $a$ and $b$.  These components are denoted $P_{a, b} (c^i)$.  See Figure \ref{barbellfig}(a).  Note that these two sets are open.  The letter $P$ is chosen here because the boundaries of these sets are taxicab parabolas.  For example, the boundary of $P_{a,b} (c^1)$ is a taxicab parabola where the focus is $c^1$ and the directrix is the guide line through $c^2$ that does not include $c^1$.

The \textit{lightning bolt of type 1} $lb^1(a, b)$ is the intersection of $bb(a, b)$ and the closed vertical strip $R_2 \cup R_5 \cup R_8$.  Similarly, The \textit{lightning bolt of type 2} $lb^2(a, b)$ is the intersection of $bb(a, b)$ and the closed horizontal strip $R_4 \cup R_5 \cup R_6$.  See Figures \ref{barbellfig}(b) and (c).  Note that the boundary of a barbell is the union of the corresponding lightning bolts.

Finally, given two points $p$ and $q$, we define the \textit{midpoint} of $p$ and $q$ to be the point
\[
m(p,q) = \left(\frac{p_1 + q_1}{2}, \frac{p_2 + q_2}{2}\right).
\]
This is the same as the midpoint of the line segment defined by $p$ and $q$ in Euclidean geometry.

\subsection{Isometries}

Recall that a map $\varphi: (X,d) \to (X,d)$ is an isometry if for all $x,y \in X$, $d(\varphi(x), \varphi(y))= d(x,y)$.  Isometries of taxicab space form a group under composition isomorphic to $\bR^2 \rtimes D_4$, \cite{isometries1}.

It follows from \cite{isometries1,isometries2} that the isometries of the taxicab plane are generated by translations by arbitrary vectors and reflections through coordinate lines and guide lines, and include rotations about any point through  Euclidean reference angles that are integer multiples of $\pi/2$, 

Finally, in the sequel we say that points $p$ and $q$ in $\bR^2$ are in standard position if $p = (0,0)$ and $q = (q_1, q_2)$ lies in the first quadrant, and $0 \leq q_2 \leq q_1$.  The following lemma will be used repeatedly:

\begin{lemma}
Given $p,q \in \mathbb{R}^2$ there exists an isometry $\varphi$ such that $\varphi(p)$ and $\varphi(q)$ are in standard position.
\end{lemma}

\begin{proof}
Consider two points $p,q \in \bR^2$.  First, apply a translation by $-p$.  Then, depending on the location of the image of $q$ under this translation, apply up to two reflections across coordinate lines or guide lines through the origin.  A choice of such reflections can always be made such that after composition with the translation, the image of $q$ lies in the first quadrant and below $gl^+(0, 0)$.
\end{proof}

\section{Apollonian sets} \label{sec:as}

Given two points $p \neq q$ in $\mathbb{R}^2$, and an extended real number number $k \in [0, \infty] $, the \textit{Apollonian set} $A(p, q; k)$ is
\[
A(p, q; k) = \left\{x \in \mathbb{R}^2: \frac{d(x,p)}{d(x,q)} = k \right\}
\]
with the convention that $A(p, q; \infty) = \{q\}$.  The case $p = q$ is avoided simply because it is somewhat degenerate:  the Apollonian sets would be empty for all $k \neq 1$ and when $k = 1$, the set would be the entire plane (with the convention that $\frac{0}{0} = 1$). 

If the Euclidean distance is used in this definition, this set is a circle if $k \in (0, 1) \cup (1, \infty)$, a line if $k = 1$, and a point if $k = 0$ or $\infty$.  This result is credited to Apollonius, see \cite{Smart}.  If $d$ is the taxicab distance, then this set can take on a number of different shapes depending on the relative positions of $p$ and $q$, and the value of $k$.  See Figure \ref{examplefig} for some typical shapes that arise.  The primary goal of this paper is to completely understand these sets.

With the exceptions of the cases where $k$ is 0, 1, or $\infty$, we will find that Apollonian sets are piecewise-linear simple closed curves, and we will need to make use of the compact sets that they bound.  With this in mind, we define the \textit{filled Apollonian set} as follows:
When $k > 1$
\begin{align*}
B(p, q; k) &= \left\{x \in \mathbb{R}^2: \frac{d(x,p)}{d(x,q)} \geq k \right\} \\
		&= \bigcup_{\kappa \in [k, \infty]} A(p, q; \kappa),
\end{align*}
and when $k < 1$
\begin{align*}
B(p, q; k) &= \left\{x \in \mathbb{R}^2: \frac{d(x,p)}{d(x,q)} \leq k \right\} \\
		&= \bigcup_{\kappa \in [0, k]} A(p, q; \kappa).
\end{align*}
Note that we do not try to define a filled Apollonian set when $k = 1$.

\subsection{Affine sets} \label{sec:affine}

More broadly, the defining relation for Apollonian sets may be cast in a slightly larger family, which in Euclidean space unifies the equations of circles, ellipses, hyperbolas, and the circle of Apollonius.  Given any  focal points $p, q \in \mathbb{R}^2$, and real parameters $\alpha, \beta, \gamma$, the affine set $S(p, q; \alpha, \beta, \gamma)$ is defined by a certain linear combination of the distances from the focal points
\[S = S(p,q;\alpha,\beta,\gamma)  = \{(x_1, x_2) \in \mathbb{R}^2 : \A \, d(x,p) + \B \, d(x,q) = \G  \}.\]
For example an Apollonian set corresponds to $\alpha = 1, \beta = -k, \gamma = 0$.  

More exploration of the affine set in taxicab space is warranted.  Without trying to perform a complete analysis, here are some general observations, true for all affine sets, that help our understanding of Apollonian sets being studied in this paper.

Written out explicitly, we see the affine set is determined by the set of points for which
\begin{equation} \label{affineeqn}
\A (|x_1 - p_1| + |x_2 - p_2|) + \B (|x_1 - q_1| + |x_2 - q_2|) = \G.
\end{equation}
The solution set can be determined by focusing on each of the 9 regions $R_i$, resolving the absolute values appropriately.  In general, when restricting to each region, the solution set is either empty, a line segment, a ray, or occasionally the entire region.

\begin{lemma} \label{filledregionlemma}
Suppose $S(p, q; \alpha, \beta, \gamma) \cap R_i = R_i$.  Then $\alpha = \beta = \gamma = 0$ or $i$ is odd.  Furthermore, if $i \neq 5$ then $\alpha = -\beta$ and if $i = 5$ then $\alpha = \beta$.
\end{lemma}

\begin{proof}
The only way every point in a region could also be in a particular affine set is if Equation \eqref{affineeqn} becomes vacuous and this can only happen if all instances of $x_1$ and $x_2$ cancel.

The odd regions are characterized by the fact that each coordinate of the points lie outside those of $p$ and $q$ or, in only the case of $R_5$, both coordinates lie between those of $p$ and $q$.  The only way then that the necessary cancelation occurs is if $\alpha = -\beta$ when $i = 1, 3, 7, 9$ or if $\alpha = \beta$ when $i = 5$.

If $i$ is even, then for any $x \in R_i$, one coordinate lies outside those of $p$ and $q$ and one coordinate lies between those of $p$ and $q$.  As a consequence, if $\alpha$ and $\beta$ are not both zero, then if $\alpha$ and $\beta$ cause one coordinate to cancel, the other coordinate must remain and so Equation \eqref{affineeqn} cannot be vacuous.  If $\alpha = \beta = 0$, then the equation is inconsistent unless $\gamma$ is also zero.
\end{proof}

It is worth pointing out that the indicated relationship between $\alpha$ and $\beta$ is not sufficient for a given odd region to be contained in an affine set.  There is also a relationship between $\alpha$, $\gamma$, $p$ and $q$ that must be satisfied.  We leave it to to the reader to determine the appropriate relationship for each region.

\begin{lemma} \label{oddregionslopelemma}
Let $i$ be odd.  If $S(p, q; \alpha, \beta, \gamma) \cap R_i$ is a line segment, then it is part of a guide line.  More specifically, for $i = 1, 9$, the slope is $+1$, for $i = 3, 7$, the slope is $-1$, and for $i = 5$, the sign of the slope is opposite that of the segment connecting $p$ and $q$.
\end{lemma}

\begin{proof}
The odd regions are characterized by the fact that each coordinate of the points lie outside those of $p$ and $q$ or, in only the case of $R_5$, both coordinates lie between those of $p$ and $q$.  This causes the absolute values to resolve in such a way that the slope of the line is always $\pm 1$.  For example, in $R_3$ the equation for $S(p, q; \alpha, \beta, \gamma)$ becomes
\[
\alpha \bigl((x_1 - p_1) + (x_2 - p_2)\bigr) + \beta \bigl((x_1 - q_1) + (x_2 - q_2)\bigr) = \gamma
\]
which reduces to
\[
x_1 + x_2 = \frac{\gamma + \alpha(p_1 + p_2) + \beta (q_1 + q_2)}{\alpha + \beta}.
\]
This is the equation of a line with slope $-1$ unless $\alpha = - \beta$, in which case the region is either empty or completely filled.
\end{proof}

\begin{lemma} \label{evenregionslopelemma}
Let $i$ be even and set $S = S(p, q; \alpha, \beta, \gamma)$.  Suppose $S \cap R_i$ is a line segment with slope $m \in [-\infty, \infty]$.  If $S \cap R_{10 - i}$ is not empty, then this set is a line segment with slope $-m$.  Also, If $j$ is even and not equal to $i$ or $10 - i$, and $S \cap R_{j}$ is not empty, then this set is a line segment with slope $\pm\frac{1}{m}$, where the sign depends on the particular region and the relative positions of $p$ and $q$.  Taken together, line segments in even numbered regions have slopes of absolute value $m$ or $1/m$.
\end{lemma}

The proof of this lemma is similar to the previous lemma.  In the even regions, the absolute values resolve to produce equations for lines with the indicated slopes, specifically of the form $\pm \frac{\alpha \pm \beta}{\alpha \mp \beta}$.  Note that if the slope is 0 or $ \pm \infty$ then the solution set may be a ray, and not just a segment.

Taken together, Lemmas \ref{filledregionlemma}, \ref{oddregionslopelemma}, and \ref{evenregionslopelemma} describe what happens when regions are not empty but determining whether or not a region is empty is subtle.  It can happen in one of two ways.  Either Equation \eqref{affineeqn} reduces to an equation that has no solution, or it reduces to a line that does not intersect the region in question.

Finally, the way in which an affine set interacts with one region can provide information about how the affine set interacts with adjacent regions.  We leave most of the cases to the reader, but make the following observation which will be useful later.

\begin{lemma} \label{sharedboundarylemma}
Let $S = S(p, q; \alpha, \beta, \gamma)$ be an affine set and let $R_i$ and $R_j$ be adjacent regions with $R_i \cap R_j = E$.  Suppose $S \cap R_i$ is a line segment that meets $E$ only at a point $b$ that is distinct from $p$, $q$, and their coordinate complements.  Then $S \cap R_j$ is also a line segment that meets $E$ at $b$.
\end{lemma}

\begin{proof}
Since $b \in E \subset R_j$, $S \cap R_j$ cannot be empty.  But $S \cap R_j$ cannot be all of $R_j$ either since that would imply that $S \cap E = E$, a contradiction.  The only remaining possibility is that $S \cap R_j$ is a line segment.  It cannot intersect $E$ at a point other than $b$ because again that would imply that $S \cap E$ consists of more than a single point.
\end{proof}

\subsection{Preliminary facts about Apollonian sets}

We now focus specifically on Apollonian sets and document some basic facts.  First, as indicated in Figure \ref{familiesfig}, choosing any $x \in \bR^2$ and computing the ratio of distances $d(x,p)/d(x,q)$ immediately yields the following lemma:

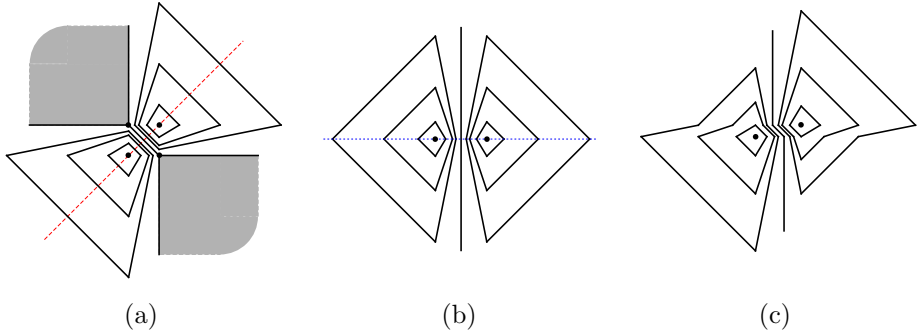
\begin{figure}
\begin{picture}(360,120)
\put(0,10){
\includegraphics[scale = .3, clip = true, draft = false]{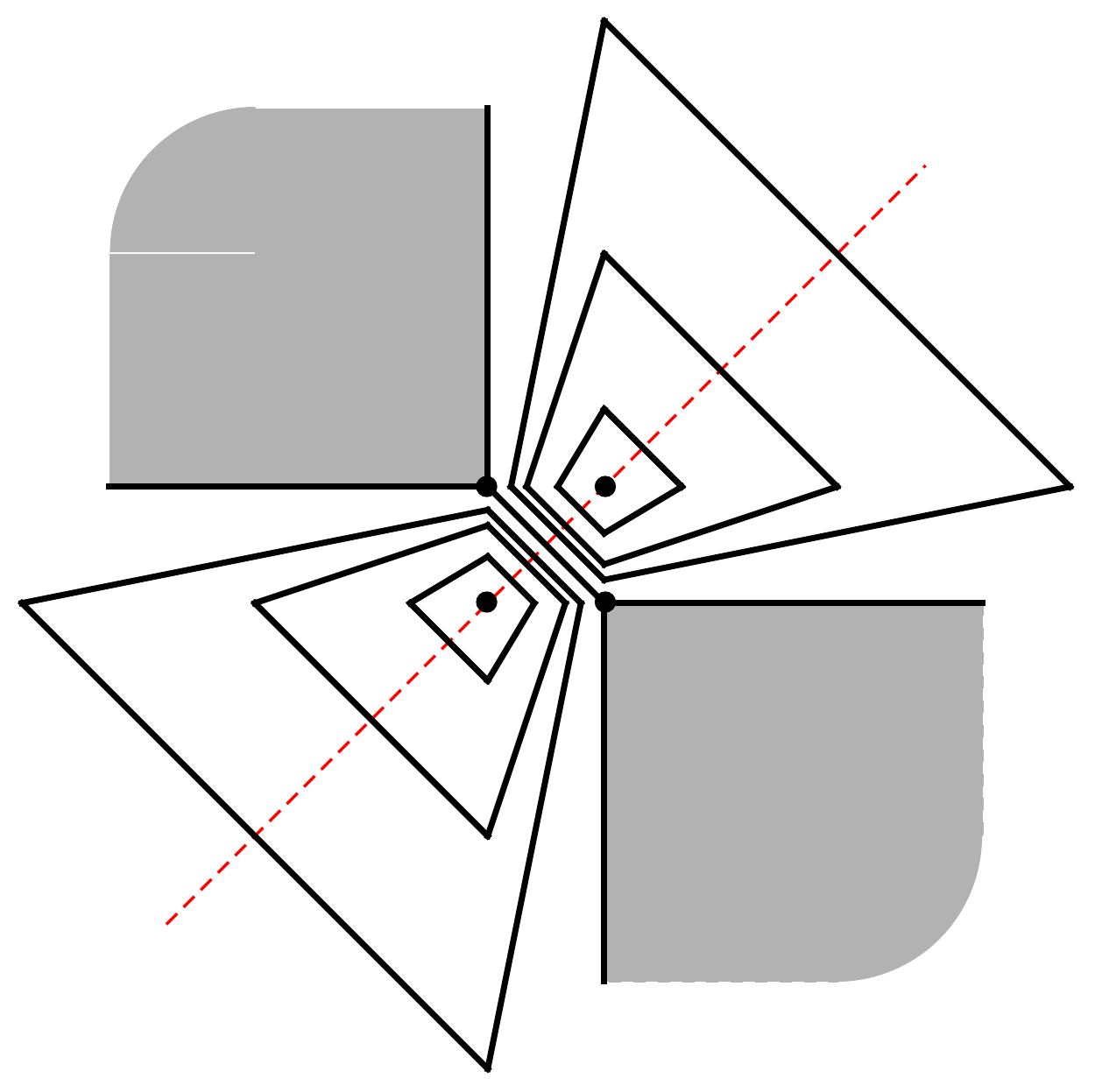}
}
\put(50,-5){(a)}

\put(120,20){
\includegraphics[scale = .3, clip = false, draft = false]{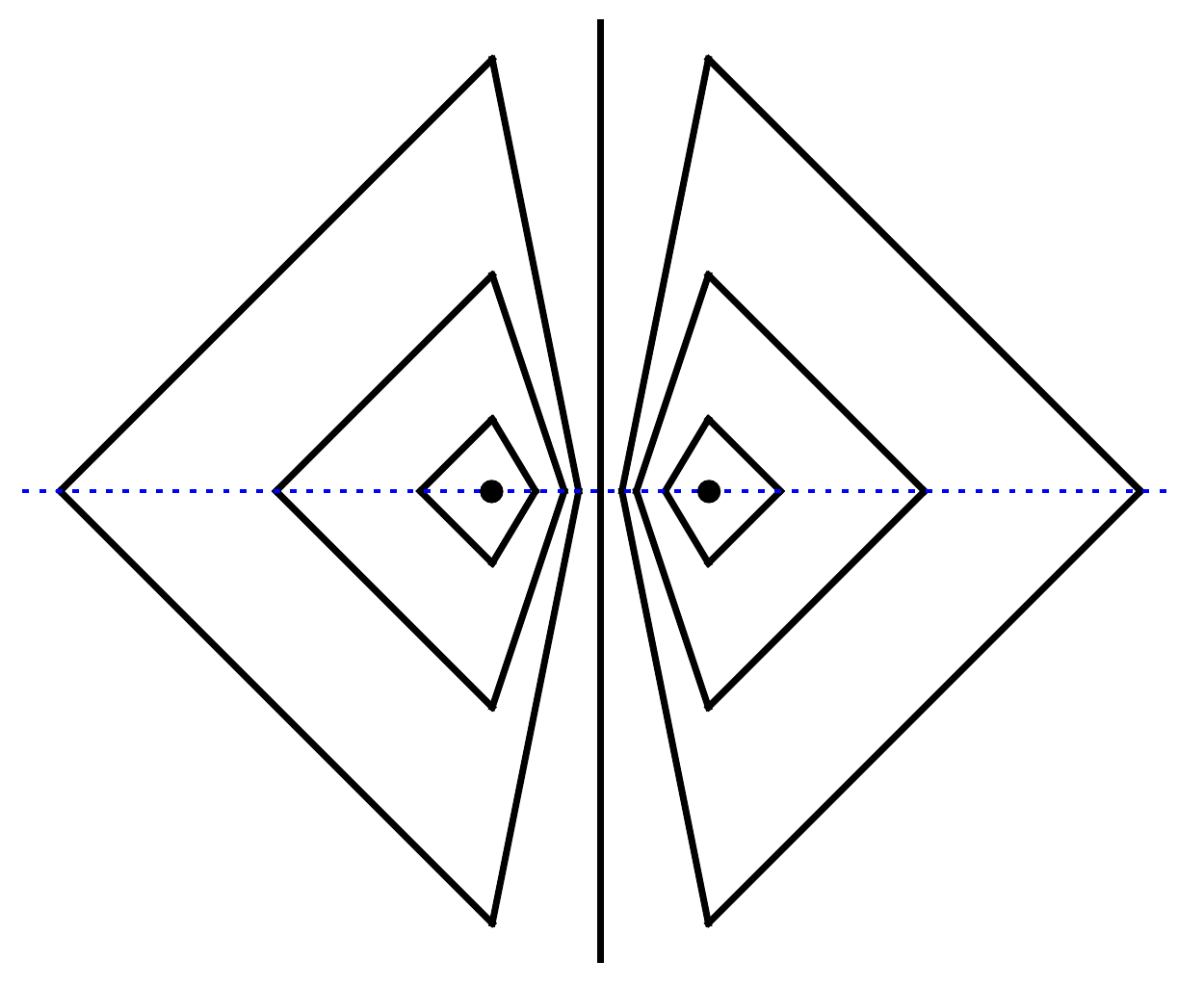}
}
\put(170,-5){(b)}

\put(240,20){
\includegraphics[scale = .3, clip = false, draft = false]{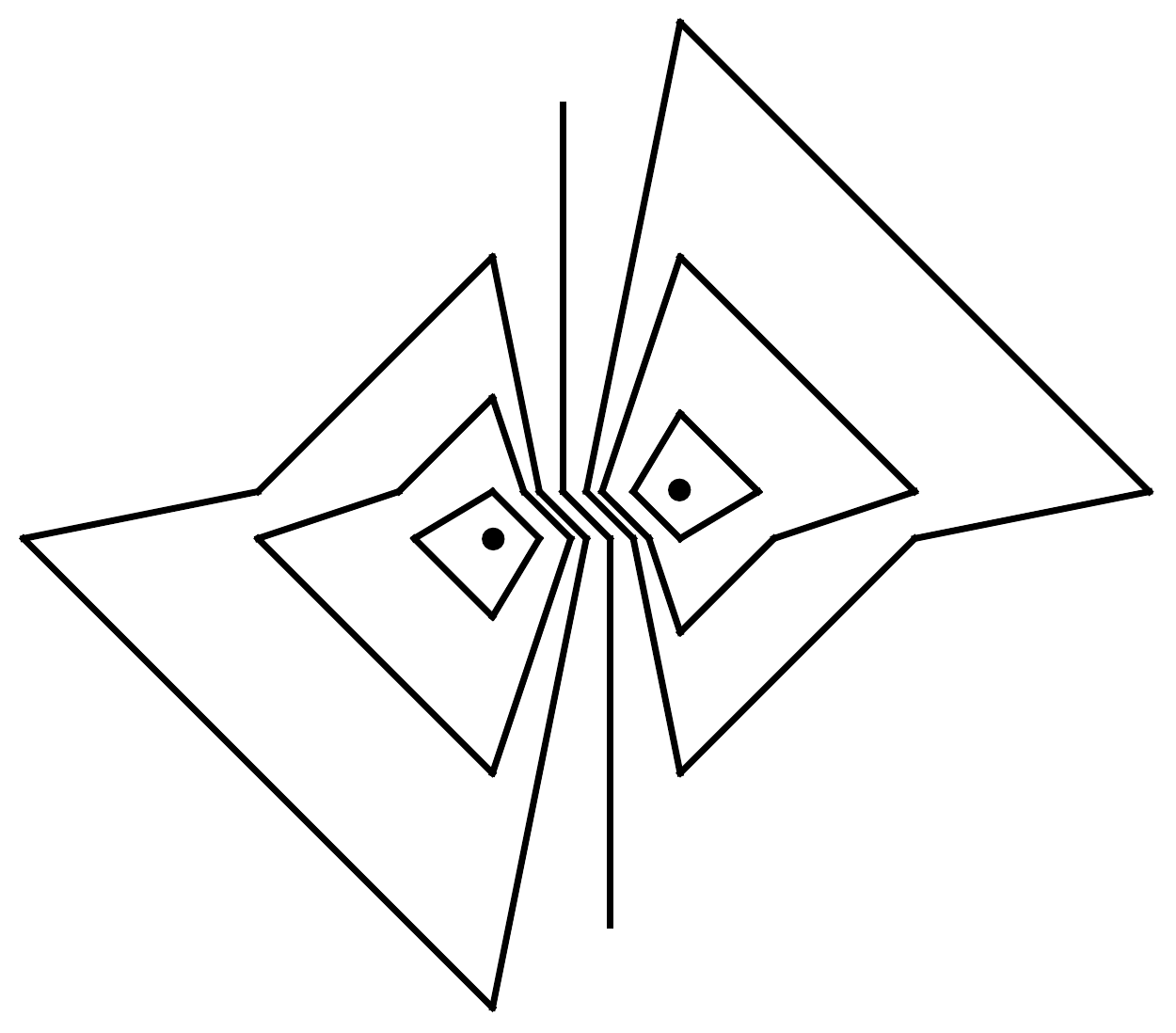}
}
\put(290,-5){(c)}

\end{picture}

\caption{Families of Apollonian sets: (a) when $p$ and $q$ share a guide line, (b) when $p$ and $q$ share a coordinate line, (c) when $p$ and $q$ share neither a guide line nor a coordinate line.  In all cases $k = 0, \frac{1}{4},\ \frac{1}{2},\ \frac{2}{3},\ 1,\ \frac{3}{2},\ 2,\ 4, \infty$ are shown.} \label{familiesfig}
\end{figure}

\begin{lemma} \label{uniquelemma}
Given $p, q \in \bR^2$, for all $x \in \bR^2$, there exists a unique $k \in [0, \infty]$ such that $x \in A(p, q; k)$.
\end{lemma}

We also note that $R_5$ is never empty.

\begin{lemma} \label{regionfivelemma}
Given $p, q \in \bR^2$, $A(p,q; k) \cap R_5$ is always a line segment.
\end{lemma}

\begin{proof}
The line segment connecting $p$ and $q$ will always have exactly one point on it satisfying the condition that defines $A(p, q; k)$.
\end{proof}.

Rearranging the formula in the definition of Apollonian set, we also have the following lemma, which allows us to restrict attention to $k \in [1, \infty]$.

\begin{lemma}\label{lemma-recip} Given $p, q \in \bR^2$, $A(p,q;\frac{1}{k}) = A(q,p;k)$, where for $k$ we use the conventions that $\frac{1}{\infty} = 0$ and $\frac{1}{0} = \infty$.
\end{lemma}

The next lemma shows that Apollonian sets are isometry invariants.

\begin{lemma} \label{lemma-isom}
Let $p,\ q \in \bR^2$, let $k \in [0, \infty]$, and let $\varphi$ be an isometry of taxicab space.  Then
\[
\varphi\bigl(A(p, q; k)\bigr) = A\bigl(\varphi(p), \varphi(q); k\bigr)
\]
and
\[
\varphi\bigl(B(p, q; k)\bigr) = B\bigl(\varphi(p), \varphi(q); k\bigr).
\]
\end{lemma}

\begin{proof}
Let $y \in A\bigl(\varphi(p), \varphi(q); k\bigr).$  Since $\varphi$ is an isometry, it is a bijection, so there exists an $x$ such that $\varphi(x) = y$ and we can substitute to get
\begin{align*}
k &= \frac{d\bigl(y, \varphi(p)\bigr)}{d\bigl(y, \varphi(q)\bigr)} \\
	&= \frac{d\bigl(\varphi(x), \varphi(p)\bigr)}{d\bigl(\varphi(x), \varphi(q)\bigr)} \\
	&= \frac{d(x, p)}{d(x, q)}. \\
\end{align*}
This implies that $x \in A(p, q; k)$, which in turn means $y \in \varphi\bigl(A(p, q; k) \bigr)$.  Therefore
\[
A\bigl(\varphi(p), \varphi(q); k\bigr) \subset \varphi\bigl(A(p, q; k) \bigr).
\]
Reversing the argument gives us the containment in the other direction.

To prove the analogous result for $B$, one may work through a similar argument as for $A$, or one may use the alternate characterization of $B$ as a union of $A$'s.

\end{proof}

\begin{lemma}\label{lemma-180}
Given $p$ and $q$, let $\varphi$ denote the isometry that rotates the plane through the angle $\pi$ about $m(p,q)$.  Then
\[
\varphi(A(p,q;k)) = A\left(p,q;\frac{1}{k}\right).
\]
\end{lemma}

\begin{proof}
By Lemma \ref{lemma-isom}, $\varphi \Bigl(A(p,q;k)\Bigr)$ is equal to $A\Bigl(\varphi(p), \varphi(q); k \Bigr)$, which in turn is equal to $A(q, p; k)$ since $\varphi$ maps $p$ to $q$ and vice versa.  Then, by Lemma \ref{lemma-recip}, this is equal to $A\left(p, q; \frac{1}{k} \right)$.
\end{proof}

\begin{figure}
\begin{picture}(200,200)
\put(0,10){
\includegraphics[scale = .5, clip = true, draft = false]{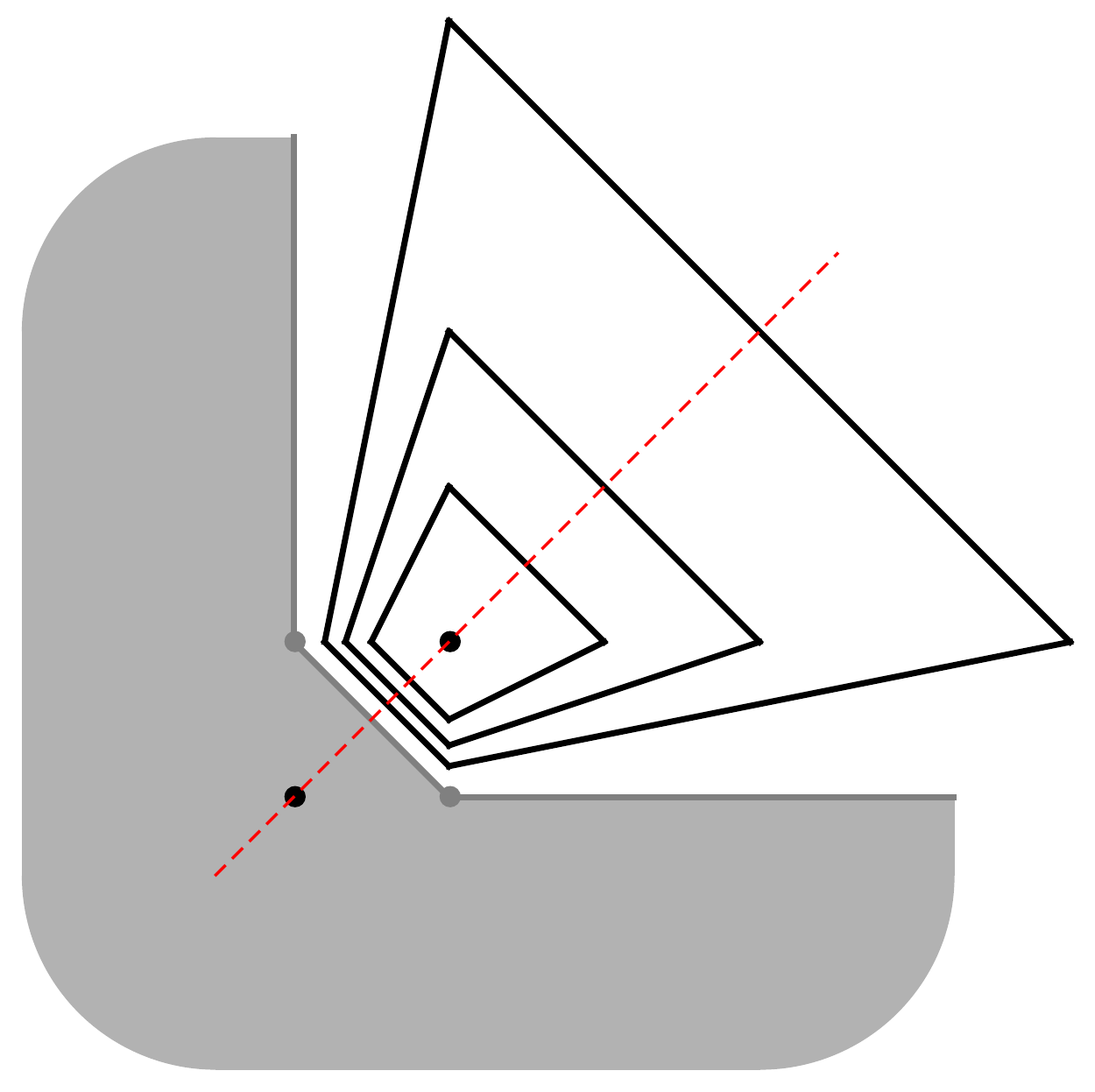}
}
\put(53,51){$p$}
\put(81,83){$q$}
\put(42,80){$c^1$}
\put(78,49){$c^2$}

\end{picture}

\caption{When $p$ and $q$ share a guide line, $A(p,q;k)$ avoids the gray area if and only $ k > 1$.  The gray area is the complement of $P_{c^1, c^2}(q)$.  The Apollonian sets shown here correspond to $k =\frac{3}{2},\ 2,\ 3$} \label{fig-filledquad}
\end{figure}
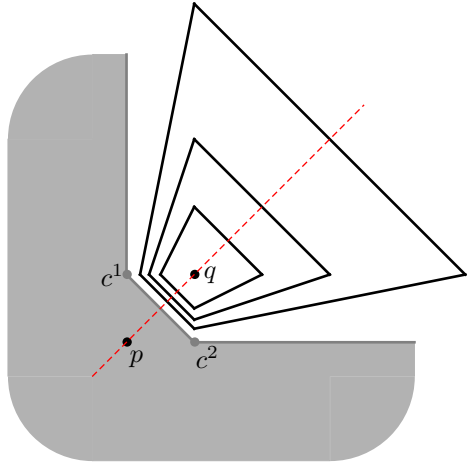

The next sequence of lemmas will provide useful bounds for the locations of Apollonian sets.  See Figure \ref{fig-filledquad}.

\begin{lemma} \label{filledquad-lemma} Let $p$ and $q$ share a guide line and let $c^1$ and $c^2$ be the coordinate complements of $p$ and $q$.  Then $x \in P_{c^1, c^2}(q)$ if and only if $\frac{d(x,p)}{d(x,q)} > 1$.
\end{lemma}

\begin{proof}
For the entirety of this proof suppose, without loss of generality, that $p$ and $q$ are in standard position.  Hence, since $p$ and $q$ share a guide line, $q_2 = q_1$.  Moreover, $P = P_{c^1,c^2}(q) = \{x \in \mathbb{R}^2: x_1 > 0,\ x_2 > 0,\ x_1 + x_2 > q_1\}$.

First, let
\[
\frac{d(x,(0,0))}{d(x,q)} >1.
\]
Thus we know 
\begin{equation} \label{loc:est1}
|x_1| + |x_2| > |x_1-q_1| + |x_2-q_1|.
\end{equation}

Our first goal is to prove that $x_1$ and $x_2$ are both positive.  We prove this by contradiction using cases.  

First suppose that both $x_1\leq0$ and $x_2\leq0$.  Then inequality \eqref{loc:est1} yields
 \[-x_1+(-x_2)>q_1-x_1+(q_1-x_2),\]
which implies $q_1 < 0$, a contradiction.

Next as $q_1=q_2$, we may without loss of generality suppose $x_1>0$ and $x_2\leq0$.  Considering all possible placements of $x_1 > 0$ relative to $q_1$ yields contradictions in inequality \eqref{loc:est1} similar to that above.

Our second goal is to prove that $x$ lies in the half plane $\{(x_1,x_2) : x_1 + x_2 > q_1\}$.  We again consider two cases, and use the fact that we now know $x_1$ and $x_2$ are both positive.  Hence, if either $x_1$ or $x_2$ is greater than $q_1$, then the result is immediate.  If both $x_1 \leq q_1$ and $x_2 \leq q_1$, then inequality \eqref{loc:est1} reduces to
\[
x_1 + x_2 > q_1 - x_1 + q_1 - x_2  
\]
and rearranging yields the result.

In the other direction, suppose $x \in P$.  Then $x$ lies in $R_2$, $R_3$, $R_5$, or $R_6$ and $x_1 + x_2 > q_1$.  It can then be checked directly that in each of these regions, $\frac{d(x,p)}{d(x,q)} > 1$.

For example, if $x \in R_5$, then $0 < x_1 \leq q_1$ and $0 < x_2 \leq q_1$ so that
\begin{equation*}
\frac{d(x,p)}{d(x,q)} = \frac{x_1 + x_2}{q_1 - x_1 + q_1 - x_2}
				=  \frac{x_1 + x_2}{2 q_1 - (x_1 + x_2)}
				> \frac{q_1}{2 q_1 - q_1}
				= 1.
\end{equation*}
where the inequality arises since $x_1 + x_2 > q_1$.

The analysis for the other three regions is easier and is left to the reader.
\end{proof}

Lemma \ref{filledquad-lemma} has the following immediate corollary:

\begin{lemma} \label{filledquadset-lemma} Let $p$ and $q$ share a guide line and let $c^1$ and $c^2$ be the coordinate complements of $p$ and $q$.  Then the Apollonian set $A(p, q; k)$ and the filled Apollonian set $B(p, q; k)$ are subsets of $P_{c^1, c^2}(q)$ if and only if $k > 1$.
\end{lemma}

The next two lemmas provide a slight recharacterization of these results that prove useful.  The proofs are almost immediate and are left to the reader.

\begin{lemma} \label{threequadrantlemma}
If $p$ and $q$ share a guide line and have coordinate complements $c^1$ and $c^2$, then $\bigcup_{k \in [0, 1]} A(p, q; k)$ is the complement of $P_{c^1, c^2}(q)$ and so contains the three $p$-quadrants that do not contain $q$.
\end{lemma}

Equivalently:

\begin{lemma} \label{threequadrantwithgreaterthan-lemma}
If $p$ and $q$ share a guide line and have coordinate complements $c^1$ and $c^2$, then $\bigcup_{k \in [1,\infty]} A(p, q; k)$ is the complement of $P_{c^1, c^2}(p)$ and so contains the three $q$-quadrants that do not contain $p$.
\end{lemma}

These results help us prove the following:

\begin{lemma} \label{R2lemma}
Let $g^+$ and $g^-$ be the guide complements of two points $p, q \in \mathbb{R}^2$.  Then
\[
\left[\bigcup_{k \in [0,1]} A(p, g^+; k) \right] \cup \left[\bigcup_{k \in [0,1]} A(p, g^-; k) \right] = \bR^2.
\]
\end{lemma}

\begin{proof}
By Lemma \ref{threequadrantlemma}, $\bigcup_{k \in [0,1]} A(p, g^+; k)$ contains all the $p$-quadrants that do not contain $g^+$ and $\bigcup_{k \in [0,1]} A(p, g^-; k)$ contains all the $p$-quadrants that do not contain $g^-$.  Their union must be all of $\bR^2$ since, by Lemma \ref{guidequadlemma}, $g^+$ and $g^-$ lie in different $p$-quadrants.
\end{proof}

\section{Characterization of Apollonian Sets} \label{sec:char}

In this section, we completely characterize Apollonian sets.  We focus first on the case $k = 1$.  Then we consider Those Apollonian sets $A(p, q; k)$ with $p$ and $q$ sharing a guide line, and we prove Theorem \ref{trapthm}.  Finally, we handle the general case, proving Theorem \ref{apollonianisunionthm}.

\begin{thm} \label{kequalsoneguidethm}
If $p$ and $q$ share a guide line, then $A(p, q; 1)$ is the barbell of $c^1(p, q)$ and $c^2(p, q)$.
\end{thm}

\begin{proof}
Note that
\[
A(p, q; 1) = \bigcup_{k \in [0, 1]} A(p, q; k) \cap \bigcup_{k \in [1, \infty]} A(p, q; k)
\]
which, by Lemmas \ref{threequadrantlemma} and \ref{threequadrantwithgreaterthan-lemma} implies
\begin{align*}
A(p, q; 1) &= \Bigl[P_{c^1, c^2}(q)\Bigr]^c \cap \Bigl[P_{c^1, c^2} (p) \Bigr]^c \\
		&= \Bigl[P_{c^1, c^2}(q) \cup P_{c^1, c^2} (p) \Bigr]^c \\
		&= bb(c^1, c^2)
\end{align*}
where the last line follows by definition.
\end{proof}

Next, we consider what happens if $p$ and $q$ do not share a guide line.

\begin{thm} \label{kequalsonenonguidethm}
If $p$ and $q$ do not share a guide line, let $gl$ be the guide line through $m(p, q)$ with the slope with the sign that is opposite that of the segment connecting $p$ and $q$, let $a$ and $b$ be the intersection of $gl$ with the boundary of $R_5$.  Then $A(p, q; 1)$ is the lightning bolt of $a$ and $b$ that intersects $R_5$ only on $gl$.
\end{thm}

This result provides an alternate characterization of the ``degenerate lines'' discussed in \cite{KAGO}.  While this result is not new, we mention it here as part of the complete characterization of Apollonian sets.  Another way to characterize which lightning bolt is correct is to note that if $p$ and $q$ do not share a guide line, then the coordinate rectangle has a long side and a short side.  The rays of the correct lightning bolt are always perpendicular to the long side.

Now that the Apollonian sets for $k = 1$ are characterized, we prove Theorem \ref{trapthm} which characterizes the Apollonian sets where $p$ and $q$ share a guide line and, without loss of generality, $k > 1$.  See Figure \ref{fig-filledquad} for examples of these sets and Figure \ref{symmetriesfig} illustrating the essential features established in the following:

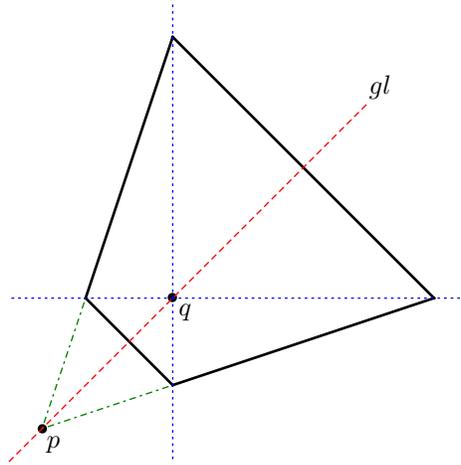
\begin{figure}
\begin{picture}(180,180)
\put(0,0){
\includegraphics[scale = .5, draft = false]{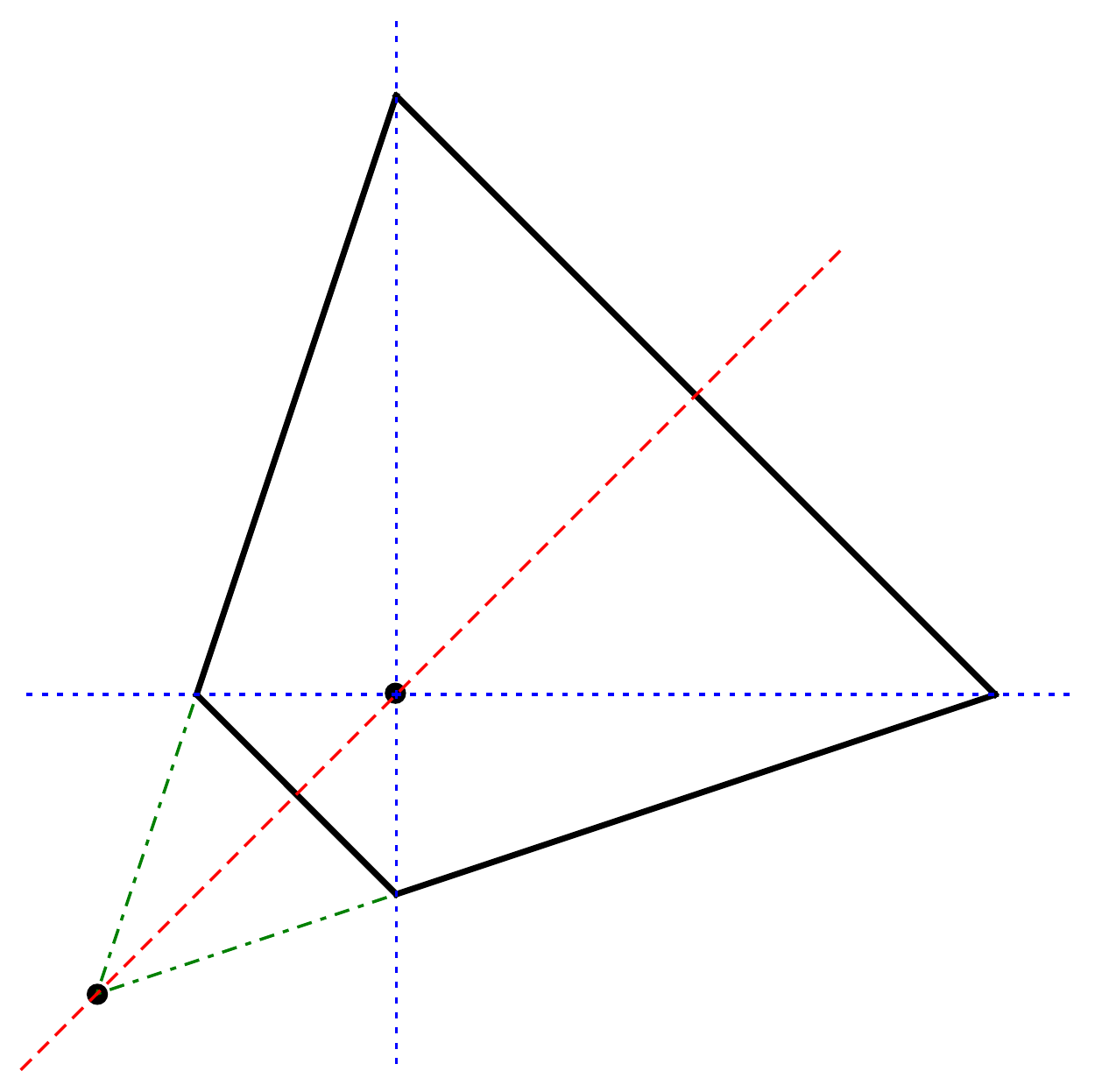}
}
\put(21,9){$p$}
\put(71,59){$q$}
\put(143,143){$gl$}

\end{picture}

\caption{The Apollonian set $A(p, q; k)$ when $p$ and $q$ share a guide line $gl$ and $k > 1$} \label{symmetriesfig}
\end{figure}

\medskip

\noindent \textbf{Theorem \ref{trapthm}.}
\emph{
Let $p$ and $q$ share a guide line $gl$ and, without loss of generality, let $k > 1$.  Then the Apollonian set $A(p, q; k)$ is an isosceles trapezoid with the following properties:
\begin{itemize}
\item the line of symmetry for $A(p, q; k)$ is $gl$;
\item the vertices of $A(p, q; k)$ all lie on the coordinate lines of $q$;
\item if the legs of $A(p, q; k)$ are extended to lines, these lines intersect at $p$;
\item the slopes of the legs are $m \frac{k+1}{k-1}$ and $m \frac{k-1}{k+1}$, where $m$ is the slope of $gl$.
\end{itemize}
}

\medskip

\begin{proof}
Without loss of generality, we may assume $p$ and $q$ are in standard position.  Then $A(p,q;k)$ is contained in the $p$-quadrant that contains $q$ by Lemma \ref{filledquad-lemma}, and thus can only have nontrivial intersection with regions $R_2, R_3, R_5$ and $R_6$.

By Lemmas \ref{oddregionslopelemma} and \ref{regionfivelemma}, $R_5$ contains a line segment with slope $-1$ which, by Lemma \ref{filledquad-lemma} must hit the boundaries shared with $R_2$ and $R_6$, implying by Lemma \ref{sharedboundarylemma} that $R_2$ and $R_6$ have nonempty intersection with $A(p, q; k)$.  Then, by Lemma \ref{evenregionslopelemma}, $R_2$ and $R_6$ contain line segments with reciprocal slopes.  By direct calculation, the equation for the segment in $R_2$ is
\[
x_2 = \frac{k+1}{k-1} x_1
\]
which has the indicated slope and passes through $p = (0,0)$.  This line will hit the boundary of $R_3$ so by  Lemmas \ref{oddregionslopelemma} and \ref{sharedboundarylemma}, the intersection of $R_3$ and $A(p, q; k)$ is a line segment with slope $-1$.

We conclude that $A(p,q;k)$ must be a quadrilateral with two parallel sides, i.e. a trapezoid. Since $p$ and $q$ are both fixed under reflection across $gl$, by Lemma \ref{lemma-isom} $\varphi( A(p,q;k) ) = A(\varphi(p), \varphi(q);k) = A(p,q;k)$ so that $A(p,q;k)$ has $gl$ as a line of symmetry.  This implies that the legs have equal length and also that when the legs are extended to lines, they both meet at $p$ since, as shown above, the leg in $R_2$ extends to a line that meets $p$.
\end{proof}

Note that, even without knowing $k$, given a single point in $A(p,q; k)$ Theorem \ref{trapthm} gives us enough information to construct the complete set.

We now come to the second main theorem of this paper, completing our characterization.

\medskip

\noindent \textbf{Theorem \ref{apollonianisunionthm}.}
\emph{
Let $p,\ q \in \bR^2$, let $g^+$ and $g^-$ be the guide complements of $p$ and $q$, and let $k \in [0, 1) \cup (1, \infty]$.  Then
\[
B(p, q; k) = B(g^+, q; k) \cup B (g^-, q; k).
\]
}

\begin{figure}
\begin{picture}(360,150)
\put(0,5){
\includegraphics[scale = .4, draft = false]{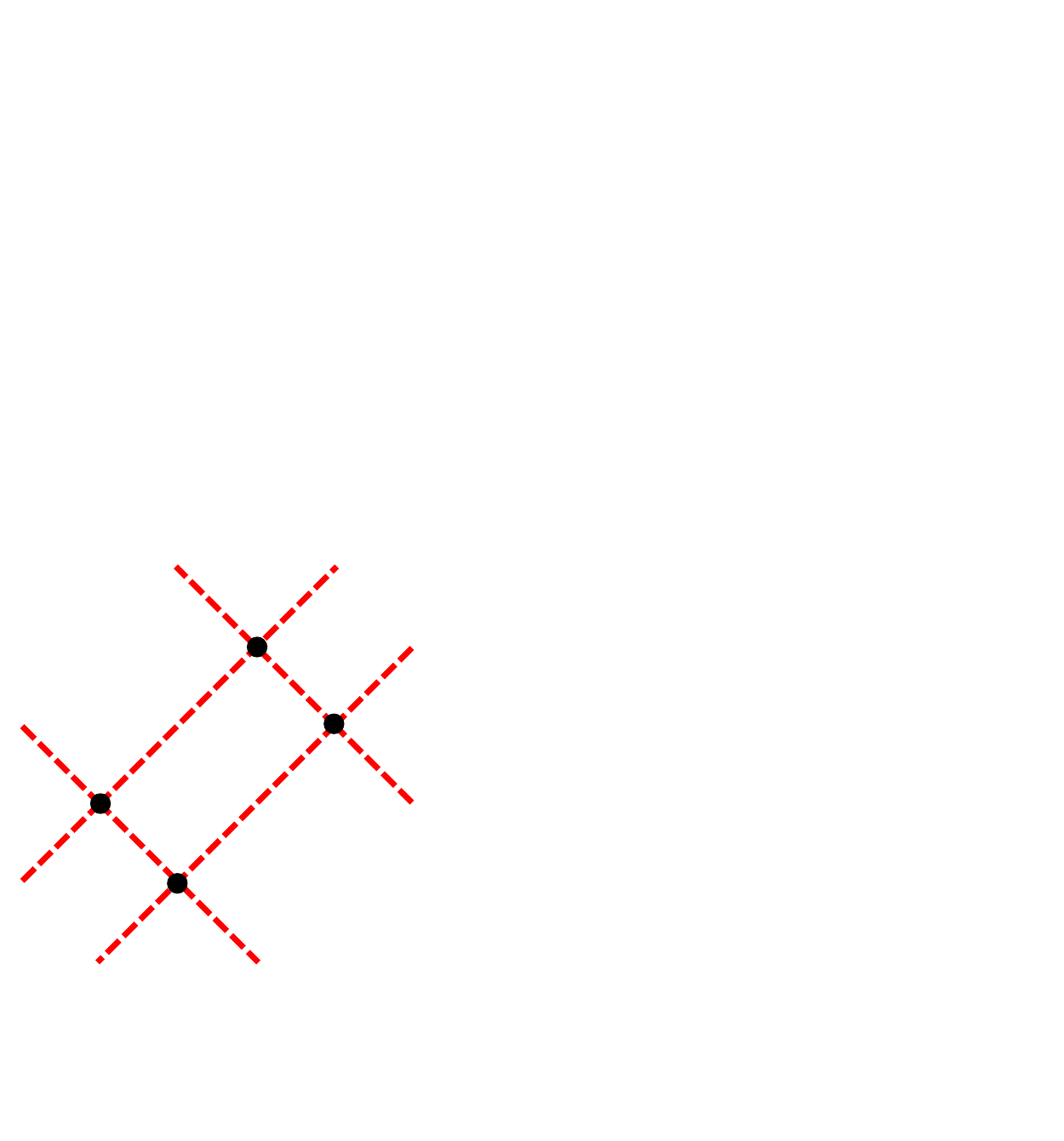}
}
\put(14,41){$p$}
\put(46,52){$q$}
\put(34,81){$g+$}
\put(22,29){$g-$}

\put(40,-5){(a)}

\put(80,5){
\includegraphics[scale = .4, draft = false]{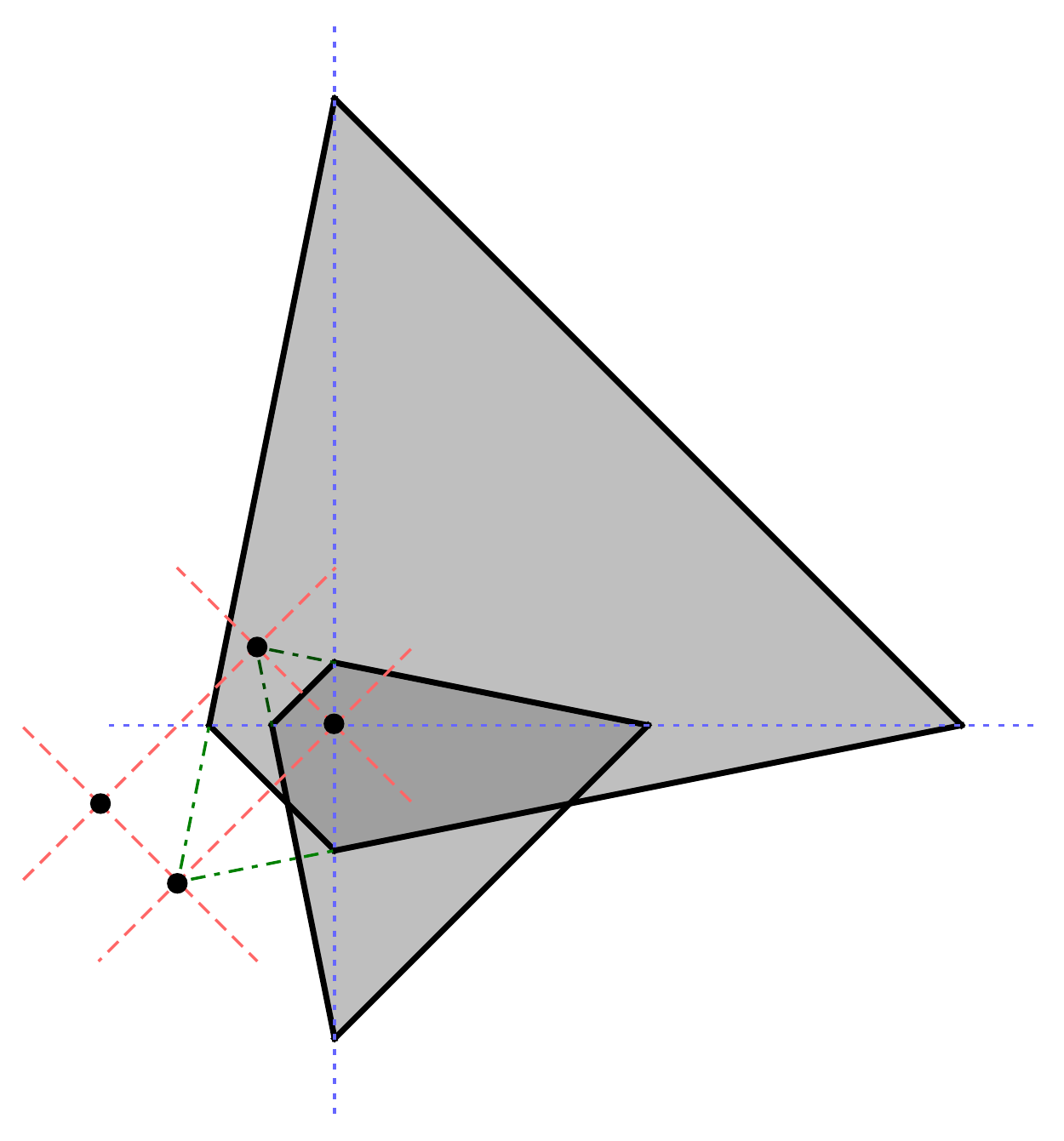}
}
\put(150,-5){(b)}

\put(220,5){
\includegraphics[scale = .4, draft = false]{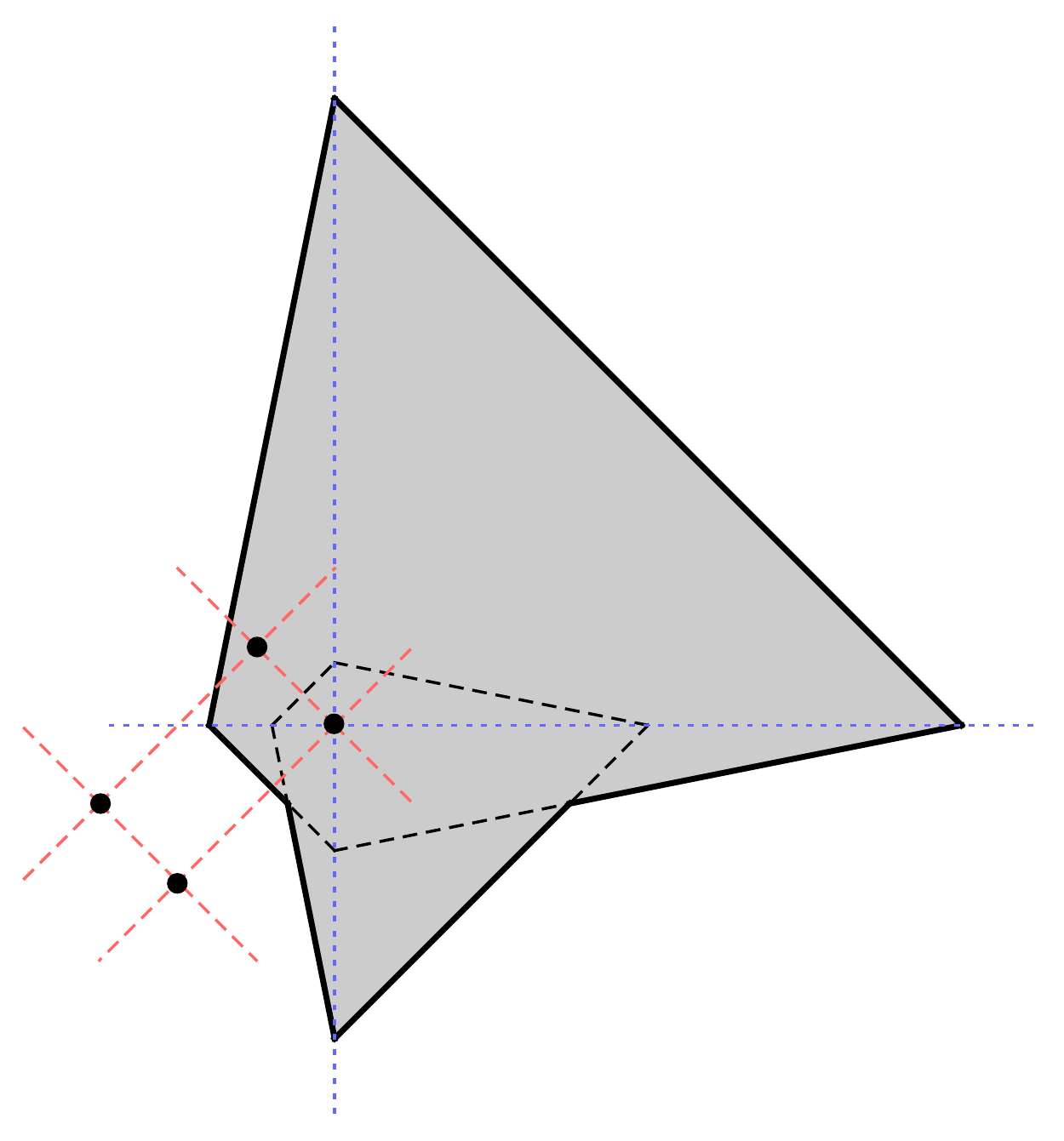}
}
\put(290,-5){(c)}

\end{picture}

\caption{Generating an Apollonian set (here, $p = (0, 0),\ q = (3, 1),\ k = \frac{3}{2}$):  (a) Find $g^+$ and $g^-$; (b) Construct the filled Apollonian sets $B(g^+, q; k)$ and $B(g^-, q; k)$; (c) The boundary of their union is $A(p, q; k)$} \label{apolloniandecompfig}
\end{figure}
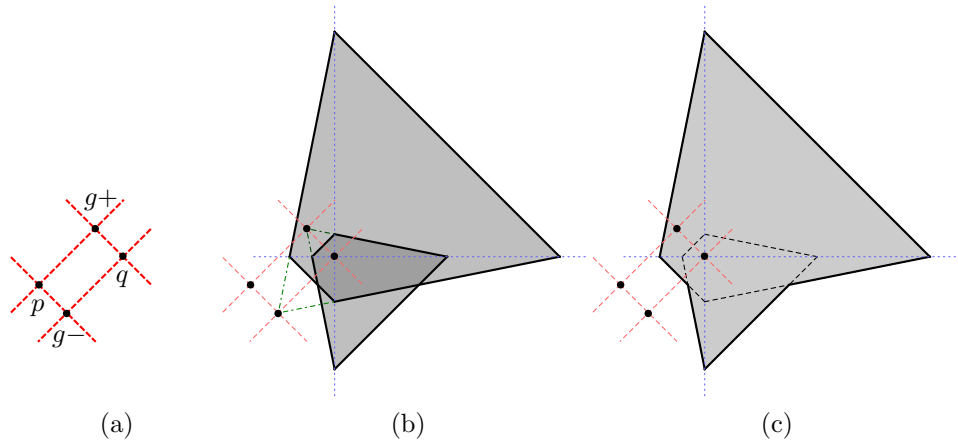

Note that the two filled Apollonian sets in the union are filled Apollonian trapezoids, so this result says that, with the exception of the case where $k = 1$, every filled Apollonian set is a union of two filled trapezoids which are themselves determined by the guide complements of $p$ and $q$.  The un-filled Apollonian set is just the boundary of this union.  See Figure \ref{apolloniandecompfig}.  As a consequence, Lemma \ref{lemma-recip} and Theorems \ref{kequalsoneguidethm}, \ref{kequalsonenonguidethm}, \ref{trapthm}, and \ref{apollonianisunionthm},  together provide a method for producing any Apollonian set in a constructive way.

\begin{proof}
Suppose, without loss of generality, that $k > 1$.  We proceed in two steps.  We first show that
\[
B(g^+, q; k) \cup B (g^-, q; k) \subset B(p, q; k).
\]
Let $x \in B(g^+, q; k)$.  Then by Lemma \ref{filledquad-lemma}, $x$ lies in the $g^+$-quadrant that contains $q$.  By Lemma \ref{guidequadlemma}, this implies that $x$ lies in one of the three $g^+$-quadrants that do not contain $p$ and so, by Lemma \ref{threequadrantwithgreaterthan-lemma},
\[
x \in \bigcup_{\kappa \in [1,\infty]} A(p, g^+; \kappa)
\]
Hence, for some $\kappa \geq 1$,
\begin{equation} \label{x-ineq-eqn}
d(x,p) \geq \kappa d(x,g^+) \geq d(x,g^+).
\end{equation}
We also know that $x \in B(g^+,q;k)$ means that
\[
\frac{d(x, g^+)}{d(x, q)} \geq k,
\]
and so $d(x, g^+) \geq k\, d(x, q)$.  Combining this with Equation \eqref{x-ineq-eqn}, we have
\[
d(x,p) \geq k\, d(x, q)
\]
and this implies $x \in B(p, q; k)$.

The argument for $x \in B(g^-, q; k)$ is exactly the same, and thus we conclude $B(g^+, q; k) \cup B (g^-, q; k) \subset B(p, q; k)$.

Next, we show that 
\[
B(p, q; k) \subset B(g^+, q; k) \cup B (g^-, q; k),
\]
by arguing the contrapositive.  If $x \not \in  B(g^+, q; k)$ and $x \not \in B(g^+, q; k)$, then we have both $d(x,g^+) < k d(x,q)$ and $d(x,g^-)< k d(x,q)$.

Now by Lemma \ref{R2lemma} applied to $p$ and $q$, there exists $\kappa \in [0,1]$ where either $x \in A(p,g^+;\kappa)$ or $x \in A(p,g^-;\kappa)$.  If $x \in A(p,g^+;\kappa)$ then
\[
d(x,p)~\leq~\kappa~d(x,g^+)~\leq~d(x,g^+),
\]
and consequently
\[
\frac{d(x,p)}{d(x,q)} = \frac{d(x,p)}{d(x,g^+)} \frac{d(x,g^+)}{d(x,q)} < k.
\]
Thus $x \not\in B(p,q;k)$, and a similar argument finishes the case when $x \in A(p,g^-;\kappa)$.
\end{proof}

\begin{figure}
\begin{picture}(290,120)
\put(0,15){
\includegraphics[scale = .35, draft = false]{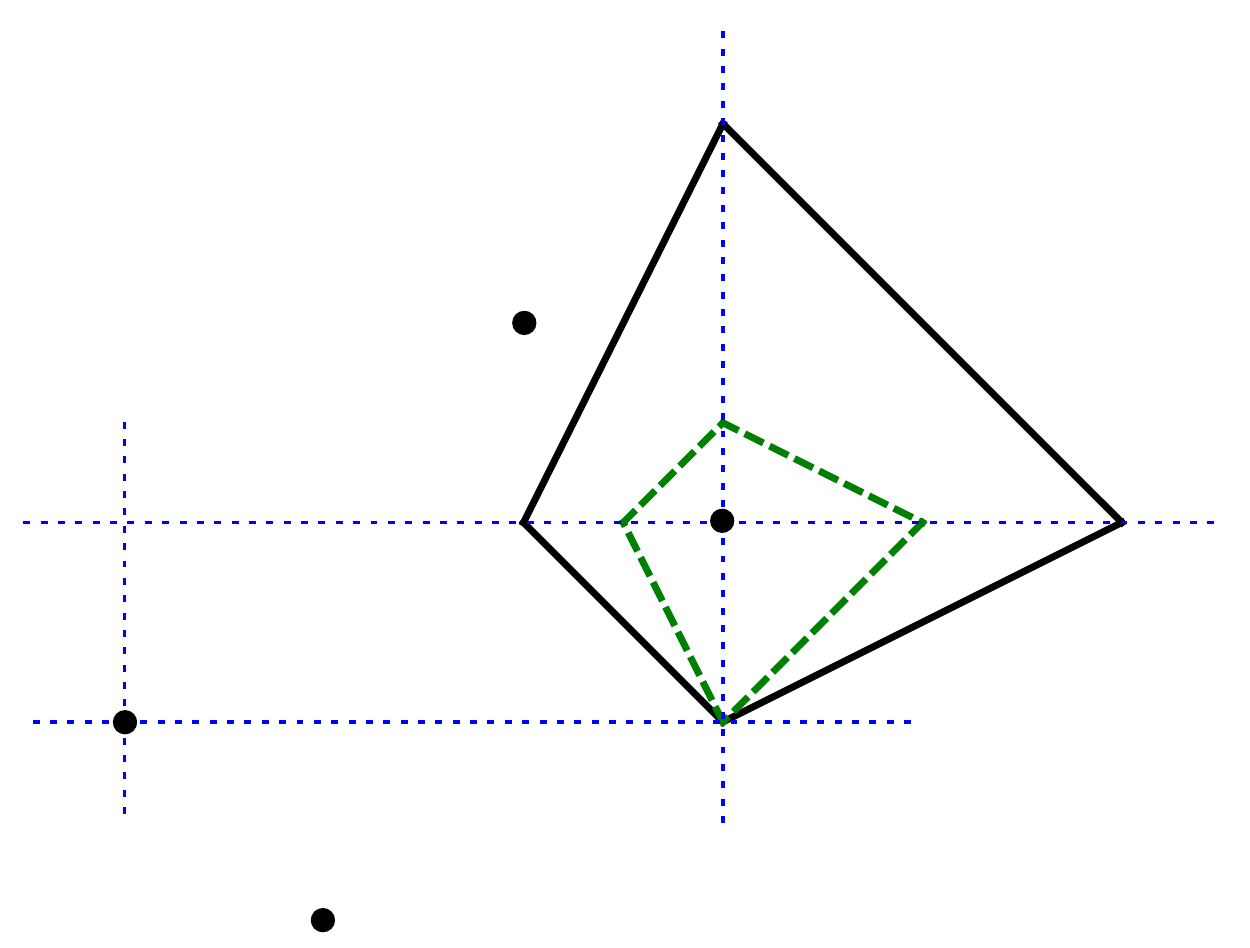}
}
\put(18,30){$p$}
\put(79,51){$q$}
\put(48,79){$g^+$}
\put(38,12){$g^-$}
\put(50,-5){(a)}

\put(160,5){
\includegraphics[scale = .35, draft = false]{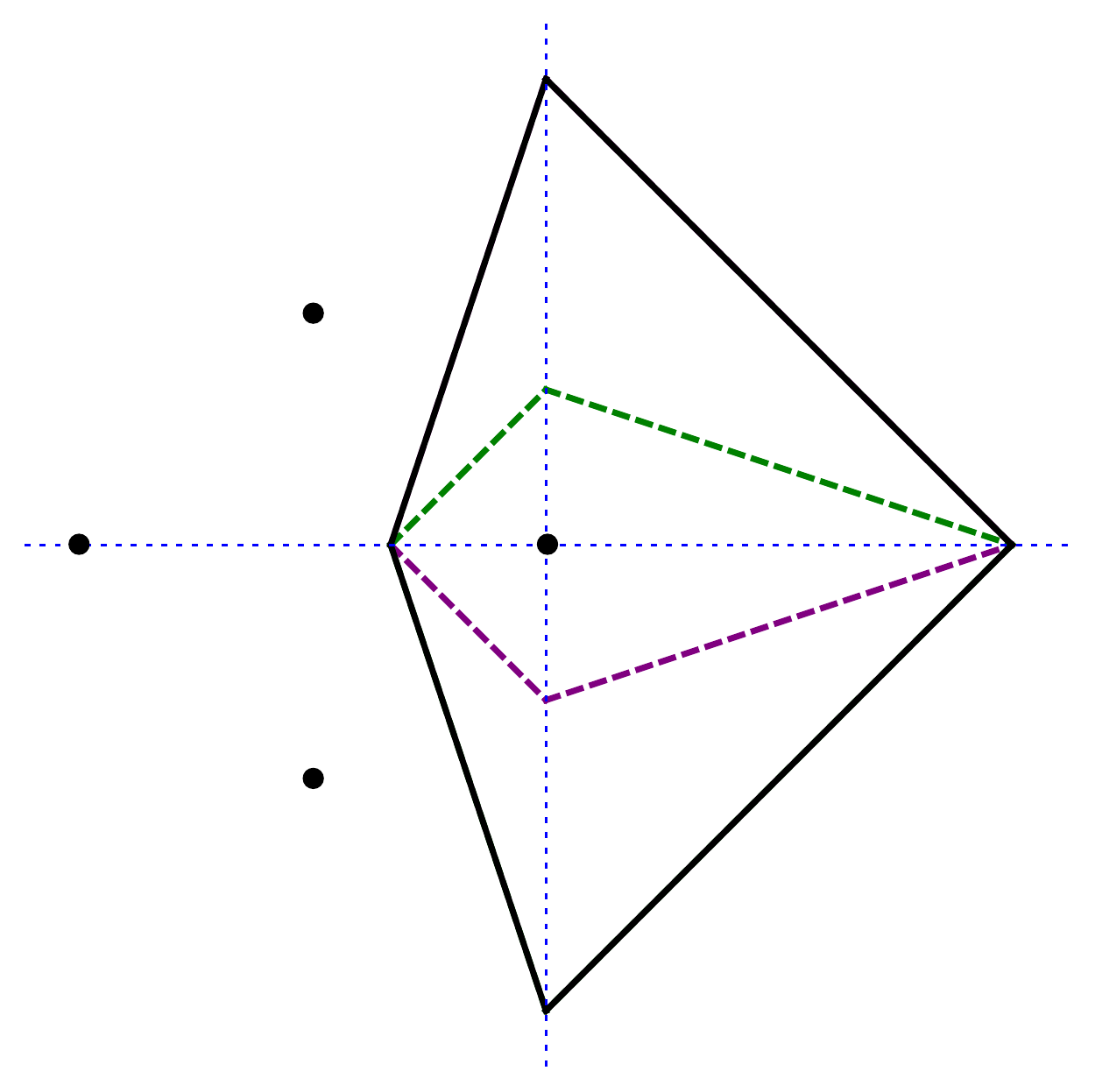}
}
\put(170,60){$p$}
\put(229,60){$q$}
\put(194,100){$g^+$}
\put(198,34){$g^-$}
\put(210,-5){(b)}

\end{picture}

\caption{Special cases: (a) The Apollonian set can be a trapezoid even when $p$ and $q$ do not share a guide line.  Here, $p = (0,0)$, $q = (3,1)$, and $k = 3$. (b) The Apollonian set is a kite when $p$ and $q$ share a coordinate line.  Here $p = (0, 0)$, $q = (4, 0)$, and $k = 2$.} \label{specialcasesfig}
\end{figure}
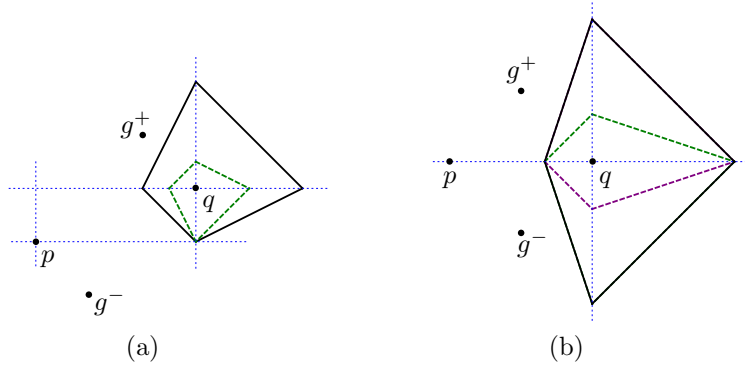

With Theorem \ref{apollonianisunionthm} established, there are a couple cases not already discussed that are of note.  First, even when $p$ and $q$ do not share a guide line, the Apollonian set can be a trapezoid as in Figure \ref{examplefig}(a).  This occurs when one of the filled Apollonian sets forming the union lies completely inside the other.  The transition between trapezoid and more complicated figure occurs when one vertex of the Apollonian set is at a coordinate complement of $p$ and $q$ as shown in Figure \ref{specialcasesfig}(a).

Second, in the special case where $p$ and $q$ share a coordinate line $cl$, the trapezoids forming the union are reflections of each other across $cl$ and the resulting figure is a kite.  See Figure \ref{specialcasesfig}(b).


\bibliographystyle{amsalpha}
\bibliography{taxicab}

\end{document}